\documentclass{amsart}

\usepackage{amsfonts,amssymb,verbatim,amsmath,amsthm,latexsym,textcomp,amscd}
\usepackage{latexsym,amsfonts,amssymb,epsfig,verbatim}
\usepackage{amsmath,amsthm,amssymb,latexsym,graphics,textcomp}
\usepackage{graphicx}
\usepackage{color}
\usepackage{url}

\input xy
\xyoption{all}

\begin{document}

\newtheorem{theorem}{Theorem}[section]
\newtheorem{prop}[theorem]{Proposition}
\newtheorem{lemma}[theorem]{Lemma}
\newtheorem{cor}[theorem]{Corollary}
\newtheorem{definition}[theorem]{Definition}
\newtheorem{conj}[theorem]{Conjecture}
\newtheorem{rmk}[theorem]{Remark}
\newtheorem{claim}[theorem]{Claim}
\newtheorem{defth}[theorem]{Definition-Theorem}

\newcommand{\boundary}{\partial}
\newcommand{\C}{{\mathbb C}}
\newcommand{\integers}{{\mathbb Z}}
\newcommand{\natls}{{\mathbb N}}
\newcommand{\ratls}{{\mathbb Q}}
\newcommand{\bbR}{{\mathbb R}}
\newcommand{\proj}{{\mathbb P}}
\newcommand{\lhp}{{\mathbb L}}
\newcommand{\tube}{{\mathbb T}}
\newcommand{\cusp}{{\mathbb P}}
\newcommand\AAA{{\mathcal A}}
\newcommand\BB{{\mathcal B}}
\newcommand\CC{{\mathcal C}}
\newcommand\DD{{\mathcal D}}
\newcommand\EE{{\mathcal E}}
\newcommand\FF{{\mathcal F}}
\newcommand\GG{{\mathcal G}}
\newcommand\HH{{\mathcal H}}
\newcommand\II{{\mathcal I}}
\newcommand\JJ{{\mathcal J}}
\newcommand\KK{{\mathcal K}}
\newcommand\LL{{\mathcal L}}
\newcommand\MM{{\mathcal M}}
\newcommand\NN{{\mathcal N}}
\newcommand\OO{{\mathcal O}}
\newcommand\PP{{\mathcal P}}
\newcommand\QQ{{\mathcal Q}}
\newcommand\RR{{\mathcal R}}
\newcommand\SSS{{\mathcal S}}
\newcommand\TT{{\mathcal T}}
\newcommand\UU{{\mathcal U}}
\newcommand\VV{{\mathcal V}}
\newcommand\WW{{\mathcal W}}
\newcommand\XX{{\mathcal X}}
\newcommand\YY{{\mathcal Y}}
\newcommand\ZZ{{\mathcal Z}}
\newcommand\CH{{\CC\HH}}
\newcommand\PEY{{\PP\EE\YY}}
\newcommand\MF{{\MM\FF}}
\newcommand\RCT{{{\mathcal R}_{CT}}}
\newcommand\PMF{{\PP\kern-2pt\MM\FF}}
\newcommand\FL{{\FF\LL}}
\newcommand\PML{{\PP\kern-2pt\MM\LL}}
\newcommand\GL{{\GG\LL}}
\newcommand\Pol{{\mathcal P}}
\newcommand\half{{\textstyle{\frac12}}}
\newcommand\Half{{\frac12}}
\newcommand\Mod{\operatorname{Mod}}
\newcommand\Area{\operatorname{Area}}
\newcommand\ep{\epsilon}
\newcommand\hhat{\widehat}
\newcommand\Proj{{\mathbf P}}
\newcommand\U{{\mathbf U}}
 \newcommand\Hyp{{\mathbf H}}
\newcommand\D{{\mathbf D}}
\newcommand\Z{{\mathbb Z}}
\newcommand\R{{\mathbb R}}
\newcommand\Q{{\mathbb Q}}
\newcommand\E{{\mathbb E}}
\newcommand\til{\widetilde}
\newcommand\length{\operatorname{length}}
\newcommand\tr{\operatorname{tr}}
\newcommand\gesim{\succ}
\newcommand\lesim{\prec}
\newcommand\simle{\lesim}
\newcommand\simge{\gesim}
\newcommand{\simmult}{\asymp}
\newcommand{\simadd}{\mathrel{\overset{\text{\tiny $+$}}{\sim}}}
\newcommand{\ssm}{\setminus}
\newcommand{\diam}{\operatorname{diam}}
\newcommand{\pair}[1]{\langle #1\rangle}
\newcommand{\T}{{\mathbf T}}
\newcommand{\inj}{\operatorname{inj}}
\newcommand{\pleat}{\operatorname{\mathbf{pleat}}}
\newcommand{\short}{\operatorname{\mathbf{short}}}
\newcommand{\vertices}{\operatorname{vert}}
\newcommand{\collar}{\operatorname{\mathbf{collar}}}
\newcommand{\bcollar}{\operatorname{\overline{\mathbf{collar}}}}
\newcommand{\I}{{\mathbf I}}
\newcommand{\tprec}{\prec_t}
\newcommand{\fprec}{\prec_f}
\newcommand{\bprec}{\prec_b}
\newcommand{\pprec}{\prec_p}
\newcommand{\ppreceq}{\preceq_p}
\newcommand{\sprec}{\prec_s}
\newcommand{\cpreceq}{\preceq_c}
\newcommand{\cprec}{\prec_c}
\newcommand{\topprec}{\prec_{\rm top}}
\newcommand{\Topprec}{\prec_{\rm TOP}}
\newcommand{\fsub}{\mathrel{\scriptstyle\searrow}}
\newcommand{\bsub}{\mathrel{\scriptstyle\swarrow}}
\newcommand{\fsubd}{\mathrel{{\scriptstyle\searrow}\kern-1ex^d\kern0.5ex}}
\newcommand{\bsubd}{\mathrel{{\scriptstyle\swarrow}\kern-1.6ex^d\kern0.8ex}}
\newcommand{\fsubeq}{\mathrel{\raise-.7ex\hbox{$\overset{\searrow}{=}$}}}
\newcommand{\bsubeq}{\mathrel{\raise-.7ex\hbox{$\overset{\swarrow}{=}$}}}
\newcommand{\tw}{\operatorname{tw}}
\newcommand{\base}{\operatorname{base}}
\newcommand{\trans}{\operatorname{trans}}
\newcommand{\rest}{|_}
\newcommand{\bbar}{\overline}
\newcommand{\UML}{\operatorname{\UU\MM\LL}}
\newcommand{\EL}{\mathcal{EL}}
\newcommand{\tsum}{\sideset{}{'}\sum}
\newcommand{\tsh}[1]{\left\{\kern-.9ex\left\{#1\right\}\kern-.9ex\right\}}
\newcommand{\Tsh}[2]{\tsh{#2}_{#1}}
\newcommand{\qeq}{\mathrel{\approx}}
\newcommand{\Qeq}[1]{\mathrel{\approx_{#1}}}
\newcommand{\qle}{\lesssim}
\newcommand{\Qle}[1]{\mathrel{\lesssim_{#1}}}
\newcommand{\simp}{\operatorname{simp}}
\newcommand{\vsucc}{\operatorname{succ}}
\newcommand{\vpred}{\operatorname{pred}}
\newcommand\fhalf[1]{\overrightarrow {#1}}
\newcommand\bhalf[1]{\overleftarrow {#1}}
\newcommand\sleft{_{\text{left}}}
\newcommand\sright{_{\text{right}}}
\newcommand\sbtop{_{\text{top}}}
\newcommand\sbot{_{\text{bot}}}
\newcommand\sll{_{\mathbf l}}
\newcommand\srr{_{\mathbf r}}
\newcommand\geod{\operatorname{\mathbf g}}
\newcommand\mtorus[1]{\boundary U(#1)}
\newcommand\A{\mathbf A}
\newcommand\Aleft[1]{\A\sleft(#1)}
\newcommand\Aright[1]{\A\sright(#1)}
\newcommand\Atop[1]{\A\sbtop(#1)}
\newcommand\Abot[1]{\A\sbot(#1)}
\newcommand\boundvert{{\boundary_{||}}}
\newcommand\storus[1]{U(#1)}
\newcommand\Momega{\omega_M}
\newcommand\nomega{\omega_\nu}
\newcommand\twist{\operatorname{tw}}
\newcommand\modl{M_\nu}
\newcommand\MT{{\mathbb T}}
\newcommand\Teich{{\mathcal T}}
\renewcommand{\Re}{\operatorname{Re}}
\renewcommand{\Im}{\operatorname{Im}}

\title{The Fourier transform on negatively curved harmonic manifolds}

\author{Kingshook Biswas}
\address{Indian Statistical Institute, Kolkata, India. Email: kingshook@isical.ac.in}

\begin{abstract} Let $X$ be a complete, simply connected harmonic manifold with sectional
curvatures $K$ satisfying $K \leq -1$, and let $\partial X$ denote the boundary at infinity
of $X$. Let $h > 0$ denote the mean curvature of horospheres in $X$, and let $\rho = h/2$.
Fixing a basepoint $o \in X$, for $\xi \in \partial X$, let $B_{\xi}$ denote the
Busemann function at $\xi$ such that $B_{\xi}(o) = 0$, then for $\lambda \in \C$
the function $e^{(i\lambda - \rho)B_{\xi}}$ is an eigenfunction of the Laplace-Beltrami operator with
eigenvalue $-(\lambda^2 + \rho^2)$.

For a function $f$ on $X$, we define the Fourier transform of $f$ by
$$\tilde{f}(\lambda, \xi) := \int_X f(x) e^{(-i\lambda - \rho)B_{\xi}(x)} dvol(x)$$ for all
$\lambda \in \C, \xi \in \partial X$ for which the integral converges. We prove a Fourier inversion
formula $$f(x) = C_0 \int_{0}^{\infty} \int_{\partial X} \tilde{f}(\lambda, \xi) e^{(i\lambda - \rho)B_{\xi}(x)} d\lambda_o(\xi) |c(\lambda)|^{-2} d\lambda$$
for $f \in C^{\infty}_c(X)$, where $c$ is a certain function on $\mathbb{R} - \{0\}$, $\lambda_o$ is the
visibility measure on $\partial X$ with respect to the basepoint $o \in X$ and $C_0 > 0$ is a constant. We also prove a Plancherel theorem.
This generalizes the corresponding results for rank one symmetric spaces of noncompact type and
negatively curved harmonic $NA$ groups (or Damek-Ricci spaces).
\end{abstract}

\bigskip

\maketitle

\tableofcontents

\section{Introduction}

\medskip

A {\it harmonic manifold} is a Riemannian manifold $X$ such that for any point $x \in X$, there exists a 
non-constant harmonic function on a punctured neighbourhood of $x$ which is radial around $x$, i.e. only depends 
on the geodesic distance from $x$. Copson and Ruse showed that this is equivalent to requiring that sufficiently small
geodesic spheres centered at $x$ have constant mean curvature, and moreover such manifolds are Einstein manifolds \cite{copsonruse}. 
Hence they have constant curvature in dimensions 2 and 3. 
The Euclidean spaces and rank one symmetric
spaces are examples of harmonic manifolds. The Lichnerowicz conjecture asserts that conversely any harmonic manifold is 
either flat or locally symmetric of rank one. The conjecture was proved for harmonic manifolds of dimension $4$ by A. G. Walker \cite{walker1}. 
In 1990 Z. I. Szabo proved the conjecture for compact simply connected harmonic manifolds \cite{szabo}.
In 1995 G. Besson, G. Courtois and S. Gallot proved the conjecture for manifolds of negative curvature admitting a
compact quotient \cite{bcg1}, using rigidity results from hyperbolic dynamics including the work of
Y. Benoist, P. Foulon and F. Labourie \cite{bfl} and that of P. Foulon and F. Labourie \cite{foulonlabourie}. In 2005 
Y. Nikolayevsky proved the conjecture for harmonic manifolds of dimension 5, showing that these must in fact 
have constant curvature \cite{nikolayevsky}. 

\medskip

In 1992 however E. Damek and F. Ricci had already provided a family of counterexamples to the conjecture in the noncompact case, which
have come to be known as {\it harmonic NA groups}, or {\it Damek-Ricci spaces} \cite{damekricci1}. These are solvable Lie groups $X = NA$ with a suitable
left-invariant Riemannian metric, given by
the semi-direct product of a nilpotent Lie group $N$ of {\it Heisenberg type} (see \cite{kaplan}) with $A = \R^+$ acting
on $N$ by anisotropic dilations. While the noncompact rank one symmetric spaces $G/K$ may be identified with harmonic $NA$ groups
(apart from the real hyperbolic spaces), there are examples of harmonic $NA$ groups which are not symmetric.
In 2006, J. Heber proved that the only complete simply connected homogeneous harmonic manifolds are the Euclidean spaces, rank one symmetric spaces,
and harmonic $NA$ groups \cite{heber}.

\medskip

Though the harmonic $NA$ groups are not symmetric in general, there is still a well developed theory of
harmonic analysis on these spaces which parallels that of the symmetric spaces $G/K$. For a noncompact symmetric space
$X = G/K$, an important role in the analysis on these spaces is played by the well-known {\it Helgason Fourier transform}
\cite{helgason1}. For harmonic $NA$ groups, F. Astengo, R. Camporesi and B. Di Blasio have defined a Fourier transform
\cite{astcampblas}, which reduces to the Helgason Fourier transform when the space is symmetric. In both cases a Fourier
inversion formula and a Plancherel theorem hold.

\medskip

The aim of the present article is to generalize these results to
general noncompact harmonic manifolds and initiate a study of harmonic analysis on these spaces.
Our analysis will apply
however only to the negatively curved harmonic manifolds. Our results may be described briefly as follows:

\medskip

Let $X$ be a complete, simply connected harmonic manifold of negative sectional curvature $K \leq -1$.
Then $X$ is a CAT(-1) space and can be compactified by adding a {\it boundary at infinity} $\partial X$,
defined in terms of equivalence classes of geodesic rays $\gamma : [0, \infty) \to X$ in $X$. Fix a basepoint $o \in X$.
Then for each boundary point $\xi \in \partial X$, there is a {\it Busemann function} $B_{\xi}$ at $\xi$ such that
$B_{\xi}(o) = 0$, defined for $x \in X$ by
$$
B_{\xi}(x) = \lim_{y \to \xi} (d(x,y) - d(o,y))
$$
where $y \in X$ converges to $\xi$ along the geodesic ray $[o, \xi)$ joining $o$ to $\xi$. The level sets of the Busemann
functions are called {\it horospheres} in $X$. The manifold $X$, being harmonic,
is also {\it asymptotically harmonic}, i.e. the mean curvature of all horospheres is equal to a constant $h \geq 0$.
Since $X$ is negatively curved, in fact $h$ is positive, and, by a result of G. Knieper \cite{knieper1}, the manifold $X$ has
{\it purely exponential volume growth}, i.e. there is a constant $C > 1$ such that for all $R > 1$ the volume
of geodesic balls $B(x, R)$ in $X$ is given by
$$
\frac{1}{C} e^{hR} \leq vol(B(x, R)) \leq C e^{hR}
$$
We let
$$
\rho = \frac{1}{2}h
$$
Then it turns out that for any $\lambda \in \C$ and $\xi \in \partial X$,
the function $f = e^{(i\lambda - \rho)B_{\xi}}$ is an eigenfunction of the Laplace-Beltrami operator $\Delta$ on $X$ with
eigenvalue $-(\lambda^2 + \rho^2)$.

\medskip

The Fourier transform of a function $f \in C^{\infty}_c(X)$ is then defined to be the function on $\C \times \partial X$ given by
$$
\tilde{f}(\lambda, \xi) = \int_X f(x) e^{(-i\lambda - \rho)B_{\xi}(x)} dvol(x)
$$
When $X$ is a noncompact rank one symmetric space, this reduces to the Helgason Fourier transform.

\medskip

The unit tangent sphere $T^1_o X$ is identified with the boundary $\partial X$ via the homeomorphism
$p : v \in T^1_o X \mapsto \xi = \gamma_v(\infty) \in \partial X$, where $\gamma_v$ is the unique geodesic
ray with $\gamma'_v(0) = v$. Pushing forward the Lebesgue measure on $T^1_o X$ by the map $p$ gives a measure
on $\partial X$ called the {\it visibility measure}, which we denote by $\lambda_o$. We then have the following
Fourier inversion formula:

\medskip

\begin{theorem} \label{mainthm1} There is a constant $C_0 > 0$ and a function $c$ on $\C - \{0\}$ such that
for any $f \in C^{\infty}_c(X)$, we have
$$
f(x) = C_0 \int_{0}^{\infty} \int_{\partial X} \tilde{f}(\lambda, \xi) e^{(i\lambda - \rho)B_{\xi}(x)} d\lambda_o(\xi) |c(\lambda)|^{-2} d\lambda
$$
for all $x \in X$.
\end{theorem}

\medskip

We also have a Plancherel formula:

\begin{theorem} \label{mainthm2} For any $f, g \in C^{\infty}_c(X)$, we have
$$
\int_X f(x) \overline{g(x)} dvol(x) = C_0 \int_{0}^{\infty} \int_{\partial X} \tilde{f}(\lambda, \xi) \overline{\tilde{g}(\lambda, \xi)} d\lambda_o(\xi)
|c(\lambda)|^{-2} d\lambda
$$
The Fourier transform extends to an isometry of $L^2(X, dvol)$ into $L^2((0, \infty) \times \partial X, C_0 d\lambda_o(\xi) |c(\lambda)|^{-2} d\lambda)$.
\end{theorem}

\medskip

The article is organized as follows. In section 2 we recall the facts about CAT(-1) spaces and harmonic manifolds
which we require. In section 3 we compute the action of the Laplacian $\Delta$ on spaces of functions constant on
geodesic spheres and horospheres respectively. In section 4 we carry out the harmonic analysis of radial functions,
i.e. functions constant on geodesic spheres centered around a given point. Unlike the well-known {\it Jacobi analysis}
\cite{koornwinder}
which applies to radial functions on rank one symmetric spaces and harmonic $NA$ groups, our analysis here is based
on {\it hypergroups} \cite{bloomheyer}. We define a spherical Fourier transform for radial functions, and obtain an
inversion formula and Plancherel theorem for this transform. In section 5 we prove the inversion formula
and Plancherel formula for the Fourier transform. The main point of the proof is an identity expressing
radial eigenfunctions in terms of an integral over the boundary $\partial X$. Finally in section 6 we define
an operation of convolution with radial functions, and show that the $L^1$ radial functions form a commutative
Banach algebra under convolution.

\medskip

{\bf Acknowledgements.} The author would like to thank Swagato K. Ray and Rudra P. Sarkar for generously
sharing their time and knowledge over the course of numerous educative and enjoyable discussions.

\medskip

\section{Preliminaries}

\medskip

In this section we recall the facts about CAT(-1) spaces and harmonic manifolds which we will require.

\medskip

\subsection{CAT(-1) spaces}

\medskip

References for the material in this section include \cite{bridsonhaefliger}, \cite{bourdon1}, \cite{bourdon2}.

\medskip

Given a real number $\kappa$, a CAT$(\kappa)$ space $X$ is a metric space with a synthetic notion
of having curvature bounded above by $\kappa$. In the following we will assume however that $\kappa \leq 0$.
It is required firstly that $X$ be a {\it geodesic metric space},
i.e. any two points $x, y$ in $X$ can be joined by a {\it geodesic}, which is an isometric embedding of an
interval of length $l = d(x,y)$ in $\mathbb{R}$ into $X$. Secondly, one requires that geodesic triangles in $X$
should be "thinner" than the corresponding triangles in the model space $M^2_{\kappa}$ of
constant curvature $\kappa$ (for $\kappa = 0$, $M^2_{\kappa} = \mathbb{R}^2$, while for $\kappa < 0$,
$M^2_{\kappa}$ is real hyperbolic space $\mathbb{H}^2$ with the distance function scaled by a factor of $1/\sqrt{-\kappa}$)
in the following sense:

\medskip

Given a geodesic triangle $\Delta$ in $X$ formed by three geodesic segments $[x,y]$, $[y,z]$, $[z,x]$, a
{\it comparison triangle} is a geodesic triangle $\overline{\Delta}$ in $M^2_{\kappa}$ formed by geodesic segments
$[\overline{x}, \overline{y}]$,  $[\overline{y}, \overline{z}]$,  $[\overline{z}, \overline{x}]$ of the same lengths as
those of $\Delta$ (such a triangle exists and is unique up to isometry). A point $\overline{p} \in [\overline{x}, \overline{y}]$
is called a {\it comparison point} for $p \in [x,y]$ if $d(x, p) = d(\overline{x}, \overline{p})$. We say $X$ is a CAT($\kappa$)
space if, in addition to being geodesic, $X$ satisfies the "CAT$(\kappa)$ inequality": for all geodesic triangles $\Delta$ in $X$
$d(p, q) \leq d(\overline{p}, \overline{q})$ for all comparison points $\overline{p}, \overline{q} \in \overline{\Delta} \subset M^2_{\kappa}$.
Alexandrov proved that a complete, simply connected Riemannian manifold is a CAT$(\kappa)$ space if and only if it
has sectional curvature bounded above by $\kappa$.

\medskip

In a CAT$(\kappa)$ space $X$, there is a unique geodesic segment $[x,y]$ joining any two points $x,y$, and moreover $X$ is contractible. Given
$x \in X$ and two points $y, z \in X$ distinct from $x$, the {\it comparison angle} $\theta_x(y, z)$ is defined to be the
angle at the vertex $\overline{x}$ in a comparison triangle in $M^2_{\kappa}$ corresponding to the geodesic triangle in $X$
with vertices $x, y, z$.

\medskip

The {\it boundary at infinity} $\partial X$ of a CAT$(\kappa)$ space $X$ is defined to be the set of
equivalence classes of geodesic rays in $X$, where a geodesic ray is an isometric embedding $\gamma : [0, \infty) \to X$
and two rays $\gamma_1, \gamma_2$ are defined to be equivalent if $\{ d(\gamma_1(t), \gamma_2(t)) : t \geq 0 \}$ is bounded
(when $X$ is a rank one symmetric space $G/K$ or a harmonic $NA$ group, then the boundary is identified with $K/M$ and
$N \cup \{\infty\}$ respectively).
For a geodesic ray $\gamma$, we denote by $\gamma(\infty) \in \partial X$ the equivalence class $[\gamma]$ of $\gamma$.
Then for any $x \in X$ and $\xi \in \partial X$ there is a unique geodesic ray $\gamma$ joining $x$ to $\xi$, i.e. such that
$\gamma(0) = x$ and $\gamma(\infty) = \xi$. We denote this ray by $[x, \xi)$.
The set $\overline{X} := X \cup \partial X$ can be given a topology called the {\it cone topology}, such that $X$ is an open dense
subset of $\overline{X}$, and for any geodesic ray $\gamma$, $\gamma(t)$ converges to $\gamma(\infty)$ as $t \to \infty$.
Moreover, $\overline{X}$ is compact if and only if $X$ is proper, i.e. closed balls in $X$ are compact.
We will assume for the rest of this section that $X$ is proper, so that $\partial X$ is compact.

\medskip

We now restrict ourselves to the case $\kappa < 0$, in which case after rescaling the metric by a constant we may assume $\kappa = -1$
so that $X$ is a CAT(-1) space.
In this case any two distinct points $\xi, \eta \in \partial X$ can be joined by a unique bi-infinite geodesic $\gamma : \mathbb{R} \to X$,
i.e. $\gamma(-\infty) = \xi, \gamma(\infty) = \eta$. We denote this bi-infinite geodesic by $(\xi, \eta)$. We remark that this
is a characteristic feature of negative curvature (or CAT(-1) spaces) as opposed to just
nonpositive curvature (or CAT(0) spaces), for example in higher rank
symmetric spaces (which are CAT(0) but not CAT(-1)) no two points in the boundary of a maximal flat can be joined by a geodesic.
For CAT(-1) spaces, we also have a natural family $\{ \rho_x : x \in X \}$ of metrics on $\partial X$ called {\it visual metrics} which are
compatible with the cone topology on $\partial X$. These are defined as follows:

\medskip

Let $x \in X$ be a basepoint. Given $\xi, \eta \in \partial X$, as points $y, z$ in $X$ converge to $\xi, \eta$
along the geodesic rays $[x, \xi), [x, \eta)$
the comparison angles $\theta_x(y, z)$ increase monotonically, and hence their limit
exists; we define the comparison angle $\theta_x(\xi, \eta)$ to be this limit.
The visual metric $\rho_x$ on $\partial X$ based at $x$ is then defined by
$$
\rho_x(\xi, \eta) = \sin\left(\frac{1}{2}\theta_x(\xi, \eta)\right)
$$
The above formula does indeed define a metric on $\partial X$, which is of diameter one. For example, if we take $X$ to be the
unit ball model $\mathbb{B}^n$ of $n$-dimensional real hyperbolic space and the point $x$ is the center of the ball,
then $\partial X$ is naturally identified with $S^{n-1}$ and the visual metric $\rho_x$ is simply half the chordal metric
on $S^{n-1}$. In general, two points $\xi, \eta \in \partial X$ are at visual distance one from each other if and only if
they are antipodal with respect to the basepoint $x$, in the sense that $x$ lies on the geodesic $(\xi, \eta)$.

\medskip

The {\it Gromov inner product} of two
points $y, z \in X$ with respect to a basepoint $x \in X$ is defined by
$$
(y | z)_x = \frac{1}{2}(d(x,y) + d(x,z) - d(y,z))
$$
When $X$ is a metric tree, the Gromov inner product equals the length of the segment common to the
geodesics $[x,y]$ and $[x,z]$. For points $\xi, \eta \in X \cup \partial X$, the Gromov inner product
$(\xi|\eta)_x$ is defined by
$$
(\xi|\eta)_x = \lim_{y \to \xi, z \to \eta} (y|z)_x
$$
(the limit exists taking values in $[0,+\infty]$, with $(\xi|\eta)_x = +\infty$ if and only if $\xi = \eta \in \partial X$).
The Gromov inner product $(\xi|\eta)_x$ measures, up to a bounded error, the distance of the geodesic $(\xi, \eta)$
from the point $x$: there exists a constant $C > 0$ such that
$$
|d(x, (\xi, \eta)) - (\xi|\eta)_x| \leq C
$$
(for $\xi \neq \eta$).
The visual metric $\rho_x$ can be expressed in terms of the Gromov inner product as
$$
\rho_x(\xi, \eta) = e^{-(\xi|\eta)_x}
$$
for $\xi, \eta \in \partial X$. It follows that for any $\epsilon > 0$, there is a constant $M = M(\epsilon) > 0$ such that if
$\rho_x(\xi, \eta) \geq \epsilon$ then the geodesic $(\xi, \eta)$ intersects the ball of radius $M$ around $x$.

\medskip

The {\it Busemann cocycle} $B : \partial X \times X \times X \to
\mathbb{R}$ is defined by
$$
B(x, y, \xi) := \lim_{z \to \xi} (d(x,z) - d(y,z))
$$
It satisfies the cocycle identity
$$
B(x, z, \xi) = B(x, y, \xi) + B(y, z, \xi)
$$
By the triangle inequality, $|B(x, y, \xi)| \leq d(x, y)$, and in fact equality occurs if and only if
$x, y$ lie on a geodesic ray with endpoint $\xi$.
For any $\xi \in \partial X$ and $x \in X$, the {\it Busemann function} at $\xi$ based at $x$
is the function on $X$ defined by $B_{\xi, x}(y) := B(y, x, \xi)$. Note that $B_{\xi, x}(x) = 0$.
The level sets of a Busemann function $B_{\xi, x}$ are called {\it horospheres} centered at $\xi$.
A key property of CAT(-1) spaces is the {\it exponential convergence of geodesics}: given two
points $p, q$ in the same horosphere centered at $\xi \in \partial X$, if $x, y$ are points on the rays
$[p, \xi), [q, \xi)$ with $d(p, x) = d(q, y) = t > 0$, then
$$
d(x, y) \leq C e^{-t}
$$
for some constant $C$ only depending on $d(p,q)$.
Thus the distance between the geodesic rays $[p, \xi), [q, \xi)$, which is a priori only bounded,
must in fact tend to zero exponentially fast.

\medskip

The visual metrics satisfy the following {\it Geometric Mean-Value Theorem}: for $x, y \in X$ and
$\xi, \eta \in \partial X$, we have
$$
\rho_y(\xi, \eta)^2 = e^{B(x,y,\xi)}e^{B(x,y,\eta)} \rho_x(\xi, \eta)^2,
$$
in particular they are all bi-Lipschitz equivalent to each other.

\medskip

We will need the following formula for Busemann functions:

\medskip

\begin{lemma} \label{buseform} For $o,x \in X, \xi \in \partial X$,
the limit of the comparison angles $\theta_o(x,z)$ exists as $z$ converges to $\xi$ along
the geodesic ray $[o,\xi)$. Denoting this limit by $\theta_o(x,\xi)$, it satisfies
$$
e^{B_{\xi,o}(x)} = \cosh(d(o,x)) - \sinh(d(o,x)) \cos(\theta_o(x,\xi))
$$
\end{lemma}

\medskip

\noindent{\bf Proof:} Consider a comparison triangle in $\mathbb{H}^2$ with side lengths $a = d(o,x), b = d(o, z), c = d(x,z)$ and angle
$\theta = \theta_0(x, z)$ at the vertex corresponding to $o$. By the hyperbolic law of cosine we have
$$
\cosh c = \cosh a \cosh b - \sinh a \sinh b \cos \theta
$$
As $z \to \xi$, we have $b, c \to \infty$, and $c - b \to B_{\xi,o}(x)$, thus
\begin{align*}
\cos \theta & = \frac{\cosh a \cosh b}{\sinh a \sinh b} - \frac{\cosh c}{\sinh a \sinh b} \\
            & \to \frac{\cosh a}{\sinh a} - \frac{e^{B_{\xi,o}(x)}}{\sinh a} \\
\end{align*}
hence the angle $\theta$ converges to a limit. Denoting this limit by $\theta_o(x, \xi)$, by the above
it satisfies
$$
e^{B_{\xi,o}(x)} = \cosh(d(o,x)) - \sinh(d(o,x)) \cos(\theta_o(x, \xi))
$$
$\diamond$

\medskip

We will need the following lemma relating the asymptotics of Busemann functions to the Gromov inner product:

\medskip

\begin{lemma} \label{busemasymptotic} Let $o \in X$.
Given $\epsilon > 0$, there exists $R = R(\epsilon) > 0$
such that for any $\xi, \eta \in \partial X$ with $\rho_o(\xi, \eta) \geq \epsilon$ and for any $x \in [o, \eta)$
with $r := d(o,x) \geq R$,
$$
|(B_{\xi,o}(x) - r) + 2(\xi|\eta)_o| \leq \epsilon
$$
In particular, for distinct points $\xi, \eta \in \partial X$,
$$
B_{\xi,o}(x) - r \to -2(\xi|\eta)_o
$$
as $x \to \eta$ along the geodesic ray $[o, \eta)$.
\end{lemma}

\medskip

\noindent{\bf Proof:} Given $\epsilon > 0$, let $\xi, \eta \in \partial X$ be such that $\rho_o(\xi, \eta) \geq \epsilon$, let
$x \in [o, \eta)$ and $r = d(o, x)$.
It is straightforward from the definitions that
$$
B_{\xi,o}(x) - r = -2(\xi|x)_o
$$
Let $x' \in X$ be the unique point of intersection of the geodesic $(\xi, \eta)$ with the
horosphere centered at $\eta$ passing through $x$, so that $B(x,x', \eta) = 0$.
Let $y, z$ be points on the rays $[o, \xi), [o, \eta)$ converging to $\xi$ and $\eta$ respectively, and
let $y', z'$ be the points of intersection of the geodesic $(\xi, \eta)$ with the horospheres based at $\xi$ and $\eta$
passing through the points $y, z$ respectively, then $d(y, y'), d(z,z') \to 0$ as $y \to \xi, z \to \eta$ by exponential convergence of geodesics.
Using $d(y', z') = d(y' , x') + d(x', z')$, it follows that $d(y, z) = d(y, x') + d(x',z) + o(1)$,
thus we have
\begin{align*}
|2(\xi|\eta)_o - 2(\xi|x)_o| & = \lim_{y \to \xi, z \to \eta} |(d(o,y) + d(o,z) - d(y, z)) - (d(o,y) + d(o, x) - d(x,y))| \\
                             & = \lim_{y \to \xi, z \to \eta} |(d(o,z) - d(o, x)) - d(y, z) + d(x,y)| \\
                             & = \lim_{y \to \xi, z \to \eta} |d(x, z) - (d(y, x')+ d(x', z) + o(1)) + d(x, y)| \\
                             & = \lim_{y \to \xi, z \to \eta} |(d(x,z) - d(x',z)) + (d(x,y) - d(x',y)) + o(1)| \\
                             & = |B(x,x',\eta) - B(x,x',\xi)| \\
                             & = |B(x,x',\xi)| \\
                             & \leq d(x, x') \\
\end{align*}
Thus it remains to estimate $d(x,x')$ for $r$ large.

\medskip

Since $\rho_o(\xi, \eta) \geq \epsilon$, there exists $M = M(\epsilon) > 0$ and a point $p \in (\xi, \eta)$ such that $d(o, p) \leq M$.
Let $p'$ be the point of intersection of the geodesic $(\xi, \eta)$ with the horosphere centered at $\eta$ passing through $o$,
so that $d(p, p') = |B(p, o, \eta)|$ and $B(p', o, \eta) = 0$.
Then $d(p', x') = d(o, x) = r$, and $d(o, p') \leq d(o, p) + d(p, p') \leq 2M$, so by exponential convergence of geodesics
there is $C  = C(M) > 0$ only depending on $M$ such that
$$
d(x, x') \leq C e^{-r}
$$
hence we can choose $R = R(\epsilon) > 0$ depending only on $\epsilon$ such that $d(x, x') \leq \epsilon$
for $r \geq R$. $\diamond$

\medskip

For $x \in X$ and $\xi \in \partial X$, let $B_x(\xi, \epsilon) \subset \partial X$ denote the ball of radius $\epsilon > 0$
around $\xi$ with respect to the visual metric $\rho_x$. The following two lemmas will be needed later in section 5 when
estimating visual measures of the visual balls $B_x(\xi, \epsilon)$.

\medskip

\begin{lemma} \label{changebasept} Let $x \in X$ and $\xi \in \partial X$. For $\epsilon > 0$, let $y$ be the point on the geodesic ray $[x, \xi)$
at a distance $r = \log(1/\epsilon)$ from $x$. Then:

\medskip

\noindent (1) For all $\eta \in B_x(\xi, \epsilon)$ we have
$$
|B(x, y, \eta) - r| \leq C
$$

\medskip

\noindent (2) There is a universal constant $\delta_0 > 0$ such that
$$
B_y(\xi, \delta_0) \subset B_x(\xi, \epsilon)
$$
\end{lemma}

\medskip

\noindent{\bf Proof:} (1): Let $\eta \in B_x(\xi, \epsilon)$. Let $z$ be the point on the geodesic ray $[x, \eta)$ at distance
$r = \log(1/\epsilon)$ from $x$, so that $d(x, y) = d(x, z) = r$. Let $\theta = \theta_x(y,z)$ be the comparison angle at $x$
between $y,z$, then $\theta \leq \theta_x(\xi, \eta)$ (by monotonicity of comparison angles along geodesics),
so $\sin(\theta/2) \leq \rho_x(\xi, \eta) \leq \epsilon$.
Applying the hyperbolic law of cosine to a comparison triangle
$\overline{\Delta}$ for the triangle $\Delta$ with vertices $x, y, z$ gives
\begin{align*}
\cosh(d(y,z)) & = 1 + 2 \sinh^2 r \sin^2\left(\frac{\theta}{2}\right) \\
              & \leq 1 + 2 (\sinh^2 r) \epsilon^2 \\
              & \leq \frac{3}{2} \\
\end{align*}
(using $\epsilon = e^{-r}$ and $(\sinh r) e^{-r} \leq 1/2$). Hence $d(y, z) \leq 2$, thus,
using $B(x, z, \eta) = r$ and the cocycle identity for the Busemann cocycle we have
\begin{align*}
|B(x, y, \eta) - r| & = |B(x, y, \eta) - B(x, z, \eta)| \\
                    & = |B(y, z, \eta)| \\
                    & \leq d(y, z) \\
                    & \leq 2 \\
\end{align*}

\medskip

(2): Let $\delta_1  = \sup \{ \delta > 0 | \overline{B_y(\xi, \delta)} \subset B_x(\xi, \epsilon) \}$. Then
$\overline{B_y(\xi, \delta_1)} \subset \overline{B_x(\xi, \epsilon)}$, and there exists
$\eta \in \overline{B_y(\xi, \delta_1)}$ such that $\rho_x(\xi, \eta) = \epsilon$. Using the Geometric Mean-Value Theorem
and part (1) above this gives
\begin{align*}
\delta^2_1 & \geq \rho_y(\xi, \eta)^2 \\
         & \geq e^{B(x, y, \xi) + B(x, y, \eta)} \rho_x(\xi, \eta)^2 \\
         & \geq e^{r + (r - 2)} e^{-2r} \\
         & = e^{-2} \\
\end{align*}
so $\delta_1 \geq e^{-1}$. Thus if we put $\delta_0 = (1/2) e^{-1}$ then $B_y(\xi, \delta_0) \subset B_x(\xi, \epsilon)$. $\diamond$

\medskip

\begin{lemma} \label{anglecomp} Let $X$ be a complete, simply connected manifold of pinched negative curvature, $-b^2 \leq K \leq -1$.
For $x \in X$ and $\xi, \eta \in \partial X$, let $\theta^R_x(\xi, \eta) \in [0, \pi]$ denote the Riemannian angle between the
geodesic rays $[x, \xi), [x, \eta)$ at the point $x$. Then
$$
\rho_x(\xi, \eta)^b \leq \sin\left(\frac{\theta^R_x(\xi, \eta)}{2}\right) \leq \rho_x(\xi, \eta)
$$
\end{lemma}

\medskip

\noindent{\bf Proof:} Since $K \leq -1$, Topogonov's theorem implies that the Riemannian angle $\theta^R_x(\xi, \eta)$ is bounded
above by the comparison angle $\theta_x(\xi, \eta)$, hence
$$
\sin\left(\frac{\theta^R_x(\xi, \eta)}{2}\right) \leq \sin\left(\frac{\theta_x(\xi, \eta)}{2}\right) = \rho_x(\xi, \eta)
$$
For the lower bound, let $y, z$ be points on the geodesic rays $[x, \xi), [x, \eta)$ respectively, and let
$\theta^b_x(y, z)$ denote the comparison angle between $y, z$ at $x$ in the model space $\mathbb{H}_{-b^2}$
of constant curvature $-b^2$, i.e. the angle in a comparison triangle in $\mathbb{H}_{-b^2}$ at the vertex corresponding to
$x$. Then the limit of the angles $\theta^b_x(y, z)$ exists as $y \to \xi, z \to \eta$, and, denoting the limit by
$\theta^b_x(\xi, \eta)$, it satisfies
$$
\sin\left(\frac{\theta^b_x(\xi, \eta)}{2}\right) = \rho_x(\xi, \eta)^b
$$
(\cite{biswas5}, Lemma 3.6). The lower curvature bound $-b^2 \leq K$ gives, by Topogonov's theorem, that the comparison
angles $\theta^b_x(y, z)$ are bounded above by the Riemannian angle $\theta^R_x(\xi, \eta)$, hence
$$
\rho_x(\xi, \eta)^b = \sin\left(\frac{\theta^b_x(\xi, \eta)}{2}\right) \leq \sin\left(\frac{\theta^R_x(\xi, \eta)}{2}\right)
$$
$\diamond$

\medskip

\subsection{Harmonic manifolds}

\medskip

References for this section include \cite{rusewalkerwillmore}, \cite{szabo}, \cite{willmore} and \cite{knieperpeyerimhoff1}.

\medskip

Throughout the rest of this paper, $(X, g)$ will denote a complete, simply connected Riemannian manifold of negative sectional curvature
$K$ satisfying $K \leq -1$, so that by the Cartan-Hadamard theorem
for any $x \in X$ the exponential map $\exp_x : T_x X \to X$ is a diffeomorphism.
Then $X$ is a CAT(-1) space with boundary $\partial X$.
It is known \cite{heintzehof} that in this case Busemann functions $B_{\xi,o}$ are $C^2$ on $X$, and $|\nabla B_{\xi, o}| = 1$,
so their level sets, the horospheres, are $C^2$ submanifolds of $X$.

\medskip

Let $T^1 X$ denote the unit tangent bundle of $X$ with fibres $T^1_x X$, $x \in X$. For any $r > 0$ we have a diffeomorphism
from the unit tangent sphere $T^1_x X$ to the geodesic sphere $S(x, r)$ of radius $r$ around $x$, $w \mapsto \exp_x(rw)$.
For any $v \in T^1_x X$ and $r > 0$, let $A(v, r)$ denote the Jacobian of this map
at the point $v$. We say that $X$ is a {\it harmonic manifold} if $A(v, r)$ does not depend on $v$, i.e. there is a
function $A$ on $(0, \infty)$ such that $A(v, r) = A(r)$ for all $v \in T^1 X$. The function $A$ is called the
{\it density function} of the harmonic manifold. The density function $A$ is increasing in $r$, and the quantity $A'(r)/A(r) \geq 0$
equals the mean curvature of geodesic spheres $S(x, r)$ of radius $r$, which decreases monotonically as $r \to \infty$ to a
constant $h \geq 0$ which equals the mean curvature of horospheres in $X$. In particular, a harmonic manifold is
{\it asymptotically harmonic}, i.e. the mean curvature of horospheres is constant.
Let $\Delta$ denote the Laplace-Beltrami operator or Laplacian on $X$,
then the mean curvature of horospheres centered at $\xi \in \partial X$ is given by $\Delta B_{\xi, o}$, so
$$
\Delta B_{\xi, o} \equiv h
$$
for all $\xi \in \partial X, o \in X$. Since Busemann functions are strictly convex in negative curvature (i.e. their Hessian
is positive definite), it follows that in fact the constant $h$ is positive.

\medskip

A function $f$ on $X$ is said to be {\it radial} around a point $x$ of $X$ if $f$ is constant on geodesic
spheres centered at $x$. For each $x \in X$, we can define a radialization operator $M_x$, defined for a
continuous function $f$ on $X$ by
$$
(M_x f)(z) = \int_{S(x, r)} f(y) d\sigma(y)
$$
where $S(x, r)$ denotes the geodesic sphere around $x$ of radius $r = d(x, z)$, and $\sigma$ denotes surface
area measure on this sphere (induced from the metric on $X$), normalized to have mass one.
The operator $M_x$ maps continuous functions to functions radial around $x$, and is formally self-adjoint,
meaning
$$
\int_X (M_x u)(z) v(z) dvol(z) = \int_X u(z) (M_x v)(z) dvol(z)
$$
for all continuous functions $u, v$ with compact support.

\medskip

 For $x \in X$, let $d_x$ denote the
distance function from the point $x$, i.e. $d_x(y) = d(x, y)$. Then $X$ is a harmonic manifold
if and only if any of the following equivalent conditions hold:

\medskip

\noindent (1) For any $x \in X$, $\Delta d_x$ is radial around $x$.

\medskip

\noindent (2) The Laplacian commutes with all the radialization operators $M_x$, i.e. $M_x \Delta u = \Delta M_x u$
for all smooth functions $u$ on $X$ and all $x \in X$.

\medskip

\noindent (3) For any smooth function $u$ and any $x \in X$, if $u$ is radial around $x$ then $\Delta u$ is radial around $x$.

\medskip

When $X$ is harmonic, note that the mean curvature $(A'/A)(r)$ of the geodesic sphere $S(x, r)$ at a point $z \in S(x, r)$ equals $\Delta d_x (z)$, hence
we have
$$
\Delta d_x = \frac{A'}{A} \circ d_x
$$

\medskip

\section{Radial and horospherical parts of the Laplacian}

\medskip

 Let $X$ be a complete, simply connected, negatively curved harmonic manifold (as before we assume sectional curvatures
 $K \leq -1$). Let $h > 0$ denote the
 mean curvature of horospheres in $X$, let $\rho = \frac{1}{2}h$, and let $A : (0, \infty) \to \R$ denote the density
 function of $X$.

\medskip

\begin{lemma} \label{laplacecomposition} For $f$ a $C^2$ function on $X$ and $u$ a $C^{\infty}$ function on
$\R$, we have
$$
\Delta(u \circ f) = (u'' \circ f) |\nabla f|^2 + (u' \circ f) \Delta f
$$
\end{lemma}

\medskip

\noindent{\bf Proof:}  Let $\gamma$ be a geodesic, then
$(u \circ f \circ \gamma)'(t) = (u' \circ f)(\gamma(t))< \nabla f, \gamma'(t) >$,
so
$$
(u \circ f \circ \gamma)''(t) = (u'' \circ f)(\gamma(t))< \nabla f, \gamma'(t) >^2 + (u' \circ f)(\gamma(t)) < \nabla_{\gamma'} \nabla f, \gamma'(t) >
$$
Now let $\{e_i\}$ be an orthonormal basis of $T_x X$, and let $\gamma_i$ be geodesics
with $\gamma'_i(0) = e_i$. Then
\begin{align*}
\Delta(u \circ f)(x) & = \sum_{i = 1}^n < \nabla_{e_i} \nabla (u \circ f), e_i > \\
                     & = \sum_{i = 1}^n (u \circ f \circ \gamma_i)''(0) \\
                     & = (u'' \circ f)(x) \sum_{i = 1}^n < \nabla f, e_i >^2 + (u' \circ f)(x) \sum_{i = 1}^n < \nabla_{e_i} \nabla f, e_i > \\
                     & = (u'' \circ f)(x) |\nabla f(x)|^2 + (u' \circ f)(x) \Delta f (x) \\
\end{align*}
 $\diamond$

\medskip

Any $C^{\infty}$ function on $X$ radial around $x \in X$ is of the form $f = u \circ d_x$ for some even
$C^{\infty}$ function $u$ on $\R$, where $d_x$ denotes the distance function from the point $x$, while
any $C^{\infty}$ function which is constant on horospheres at $\xi \in \partial X$ is of the form
$f = u \circ B_{\xi, x}$ for some $C^{\infty}$ function $u$ on $\R$. The following proposition says that
the Laplacian $\Delta$ leaves invariant these spaces of functions, and describes the action of the Laplacian
on these spaces:

\medskip

\begin{prop} \label{radialhorospherical} Let $x \in X, \xi \in \partial X$.

\medskip

\noindent (1) For $u$ a $C^{\infty}$ function on $(0, \infty)$,
$$
\Delta (u \circ d_x) = (L_R u) \circ d_x
$$
where $L_R$ is the differential operator on $(0, \infty)$ defined by
$$
L_R = \frac{d^2}{dr^2} + \frac{A'(r)}{A(r)} \frac{d}{dr}
$$

\noindent (2) For $u$ a $C^{\infty}$ function on $\R$,
$$
\Delta (u \circ B_{\xi, x}) = (L_H u) \circ B_{\xi, x}
$$
where $L_H$ is the differential operator on $\R$ defined by
$$
L_H = \frac{d^2}{dt^2} + 2\rho \frac{d}{dt}
$$
\end{prop}

\medskip

\noindent{\bf Proof:} Noting that $|\nabla d_x| = 1, |\nabla B_{\xi, x}| = 1$, and
$\Delta d_x = (A' / A) \circ d_x, \Delta B_{\xi,x} = 2\rho$, the Proposition follows
immediately from the previous Lemma. $\diamond$

\medskip

Accordingly, we call the differential operators $L_R$ and $L_H$
the {\it radial and horospherical parts of the Laplacian} respectively. It follows from the above
proposition that a function $f = u \circ d_x$ radial around $x$ is an eigenfunction of $\Delta$
with eigenvalue $\sigma$ if and only if $u$ is an eigenfunction of $L_R$ with eigenvalue $\sigma$.
Similarly, a function $f = u \circ B_{\xi,x}$ constant on horospheres at $\xi$ is an eigenfunction of
$\Delta$ with eigenvalue $\sigma$ if and only if $u$ is an eigenfunction of $L_H$ with eigenvalue $\sigma$.
In particular, we have the following:

\medskip

\begin{prop} \label{helgasonkernel} Let $\xi \in \partial X, x \in X$. Then for any $\lambda \in \C$, the function
$$
f = e^{(i\lambda - \rho)B_{\xi,x}}
$$
is an eigenfunction of the Laplacian with eigenvalue $-(\lambda^2 + \rho^2)$ satisfying $f(x) = 1$.
\end{prop}

\medskip

\noindent{\bf Proof:} This follows from the fact that
the function $u(t) = e^{(i\lambda - \rho)t}$ on $\R$ is an eigenfunction of $L_H$ with eigenvalue $-(\lambda^2 + \rho^2)$,
and $B_{\xi,x}(x) = 0$ gives $f(x) = 1$. $\diamond$

\medskip

\section{Analysis of radial functions}

\medskip

As we saw in the previous section, finding radial eigenfunctions of the Laplacian
amounts to finding eigenfunctions of its radial part $L_R$. When $X$ is a rank one
symmetric space $G/K$, or more generally a harmonic $NA$ group, then the volume density
function is of the form $A(r) = C \left(\sinh\left(\frac{r}{2}\right)\right)^p \left(\cosh\left(\frac{r}{2}\right)\right)^q$,
for a constant $C > 0$ and integers $p, q \geq 0$, and so the radial part $L_R = \frac{d^2}{dr^2} + (A'/A) \frac{d}{dr}$
falls into the general class of {\it Jacobi operators}
$$
L_{\alpha, \beta} = \frac{d^2}{dr^2} + ((2\alpha + 1)\coth r + (2\beta + 1)\tanh r)\frac{d}{dr}
$$
for which there is a detailed and well known harmonic analysis in terms of eigenfunctions (called {\it Jacobi functions})
\cite{koornwinder}.
For a general harmonic manifold $X$, the explicit form of the density function $A$ is not known, so it is unclear whether
the radial part $L_R$ is a Jacobi operator.
However, there is a harmonic analysis, based on hypergroups
(\cite{chebli1}, \cite{chebli2}, \cite{trimeche1}, \cite{trimeche2}, \cite{trimeche3}, \cite{bloomxu1}, \cite{xu1}),
for more general second-order differential operators on $(0, \infty)$ of the form
\begin{equation} \label{formofl}
L = \frac{d^2}{dr^2} + \frac{A'(r)}{A(r)} \frac{d}{dr}
\end{equation}
where $A$ is a function on $[0,\infty)$ satisfying certain hypotheses which allow one to endow $[0,\infty)$
with a hypergroup structure, called a {\it Chebli-Trimeche hypergroup}. We first recall some basic facts about Chebli-Trimeche hypergroups,
and then show that the density function of a harmonic manifold satisfies the hypotheses required in order to apply this theory.

\medskip

\subsection{Chebli-Trimeche hypergroups}

\medskip

A hypergroup $(K, *)$ is a locally compact Hausdorff space $K$ such that
the space $M^b(K)$ of finite Borel measures on $K$ is endowed with a product $(\mu, \nu) \mapsto \mu * \nu$
turning it into an algebra with unit, and $K$ is endowed with an involutive homeomorphism $x \in K \mapsto \tilde{x} \in K$,
such that the product and the involution satisfy certain natural properties (see \cite{bloomheyer} Chapter 1 for the
precise definition). A motivating example relevant to the following is the algebra of finite radial measures on a noncompact
rank one symmetric space $G/K$ under convolution; as radial measures can be viewed as measures on $[0, \infty)$, this endows
$[0, \infty)$ with a hypergroup structure (with the involution being the identity).
It turns out that this hypergroup structure on $[0, \infty)$ is a special case
of a general class of hypergroup structures on $[0, \infty)$ called {\it Sturm-Liouville hypergroups} (see \cite{bloomheyer}, section 3.5).
These hypergroups arise from Sturm-Liouville boundary problems on $(0, \infty)$. We will be interested in a particular class
of Sturm-Liouville hypergroups called {\it Chebli-Trimeche hypergroups}. These arise as follows (we refer to \cite{bloomheyer}
for proofs of statements below):

\medskip

A {\it Chebli-Trimeche function} is a continuous function $A$ on $[0, \infty)$ which is $C^{\infty}$ and
positive on $(0, \infty)$ and satisfies the
following conditions:

\medskip

\noindent (H1) $A$ is increasing, and $A(r) \to +\infty$ as $r \to +\infty$.

\medskip

\noindent (H2) $A'/A$ is decreasing, and $\rho = \frac{1}{2} \lim_{r \to \infty} A'(r)/A(r) > 0$.

\medskip

\noindent (H3) For $r > 0$, $A(r) = r^{2\alpha + 1} B(r)$ for some $\alpha > -1/2$ and some even, $C^{\infty}$ function $B$ on $\R$
such that $B(0) > 0$.

\medskip


Let $L$ be the differential operator on $C^2(0, \infty)$ defined by equation (\ref{formofl}),
where $A$ satisfies conditions (H1)-(H3) above. Define the differential operator $l$ on $C^2((0,\infty)^2)$
by
\begin{align*}
l[u](x,y) & = (L)_x u(x, y) - (L)_y u(x,y) \\
          & = \left(u_{xx}(x,y) + \frac{A'(x)}{A(x)}u_x(x,y)\right) - \left(u_{yy}(x,y) + \frac{A'(y)}{A(y)}u_y(x,y)\right) \\
\end{align*}

For $f \in C^2([0, \infty)^2)$ denote by $u_f$ the solution of the hyperbolic Cauchy problem
\begin{align*}
l[u_f] & = 0, \\
u_f(x, 0) = u_f(0, x) & = f(x), \\
(u_f)_y(x, 0) & = 0, \\
(u_f)_x(0, y) & = 0 \ \hbox{ for } x, y \in [0, \infty) \\
\end{align*}

For $x \in [0, \infty)$, let $\epsilon_x$ denote the Dirac measure of mass one at $x$. Then for all $x, y \in [0, \infty)$, there
exists a probability measure on $[0, \infty)$ denoted by $\epsilon_x * \epsilon_y$ such that
$$
\int_{0}^{\infty} f d(\epsilon_x * \epsilon_y) = u_f(x, y)
$$
for all even, $C^{\infty}$ functions $f$ on $\R$.  We have $\epsilon_x * \epsilon_y = \epsilon_y * \epsilon_x$ for all $x,y$, and
the product
$(\epsilon_x, \epsilon_y) \mapsto \epsilon_x * \epsilon_y$ extends to a product on all finite measures on
$[0, \infty)$ which turns $[0, \infty)$ into a commutative hypergroup $([0,\infty), *)$ (with the involution being the identity),
called the Chebli-Trimeche hypergroup associated to the function $A$.
Any hypergroup has a Haar measure, which in this case
is given by the measure $A(r) dr$ on
$[0, \infty)$.

\medskip

For a commutative hypergroup $K$ with a Haar measure $dk$, a Fourier analysis can be carried out analogous to the Fourier analysis
on locally compact abelian groups. There is a dual space $\hat{K}$ of characters, which are bounded multiplicative functions
on the hypergroup $\chi : K \to \C$ satisfying $\chi(\tilde{x}) = \overline{\chi(x)}$, where multiplicative means that
$$
\int_K \chi d(\epsilon_x * \epsilon_y) = \chi(x) \chi(y)
$$
for all $x, y \in K$. For $f \in L^1(K)$, the Fourier transform of $f$ is the function $\hat{f}$ on $\hat{K}$ defined by
$$
\hat{f}(\chi) = \int_K f \overline{\chi} dk
$$
The Levitan-Plancherel Theorem states that there is a measure $d\chi$ on $\hat{K}$
called the Plancherel measure, such that the mapping $f \mapsto \hat{f}$ extends from
$L^1(K) \cap L^2(K)$ to an isometry from $L^2(K)$ onto $L^2(\hat{K})$. The inverse Fourier transform of a function
$\sigma \in L^1(\hat{K})$ is the function $\check{\sigma}$ on $K$ defined by
$$
\check{\sigma}(k) = \int_{\hat{K}} \sigma(\chi) \chi(k) d\chi
$$
The Fourier inversion theorem then states that if $f \in L^1(K) \cap C(K)$ is such that $\hat{f} \in L^1(\hat{K})$,
then $f = (\hat{f})\check{}$, i.e.
$$
f(x) = \int_{\hat{K}} \hat{f}(\chi) \chi(x) d\chi
$$
for all $x \in K$.

\medskip

For the Chebli-Trimeche hypergroup, it turns out that the multiplicative functions on the hypergroup are
given precisely by eigenfunctions of the operator $L$.
For any $\lambda \in \C$, the equation
\begin{equation} \label{evaleqn}
Lu = -(\lambda^2 + \rho^2)u
\end{equation}
has a unique solution $\phi_{\lambda}$ on $(0, \infty)$ which extends continuously to $0$ and satisfies
$\phi_{\lambda}(0) = 1$ (note that the coefficient $A'/A$ of the operator $L$ is singular at $r = 0$ so existence of a
solution continuous at $0$ is not immediate). The function $\phi_{\lambda}$ extends to a $C^{\infty}$ even function on $\R$.
Since equation (\ref{evaleqn}) reads the same for $\lambda$ and $-\lambda$, by uniqueness we have $\phi_{\lambda} = \phi_{-\lambda}$.

\medskip

The multiplicative functions on $[0, \infty)$ are then exactly the functions $\phi_{\lambda}, \lambda \in \C$. The functions
$\phi_{\lambda}$ are bounded if and only if $| \Im \lambda | \leq \rho$. Furthermore, the involution on the hypergroup being the
identity, the characters of the hypergroup are real-valued, which occurs for $\phi_{\lambda}$ if and only if $\lambda \in \R \cup i\R$.
Thus the dual space of the hypergroup is given by
$$
\hat{K} = \{ \phi_{\lambda} | \lambda \in [0, \infty) \cup [0, i\rho] \}
$$
which we identify with the set $\Sigma = [0, \infty) \cup [0, i\rho] \subset \C$.

\medskip

The hypergroup Fourier transform of a function $f \in L^1([0, \infty), A(r) dr)$ is given by
$$
\hat{f}(\lambda) = \int_{0}^{\infty} f(r) \phi_{\lambda}(r) A(r) dr
$$
for $\lambda \in \Sigma$ (when the hypergroup arises from convolution of
radial measures on a rank one symmetric space $G/K$, then this is the well-known Jacobi transform \cite{koornwinder}).
The Levitan-Plancherel and Fourier inversion theorems for the hypergroup give the existence
of a Plancherel measure $\sigma$ on $\Sigma$ such that the Fourier transform defines an isometry from $L^2([0, \infty), A(r) dr)$ onto
$L^2(\Sigma, \sigma)$, and, for any function $f \in L^1([0, \infty), A(r) dr) \cap C([0, \infty))$ such that $\hat{f} \in L^1(\Sigma, \sigma)$,
we have
$$
f(r) = \int_{\Sigma} \hat{f}(\lambda) \phi_{\lambda}(r) d\sigma(\lambda)
$$
for all $r \in [0, \infty)$.

\medskip

In \cite{bloomxu1}, it is shown that under certain extra conditions on the function $A$, the support of the Plancherel
measure is $[0, \infty)$ and the Plancherel measure is absolutely continuous with respect to Lebesgue measure $d\lambda$ on $[0,\infty)$,
given by
$$
d\sigma(\lambda) = C_0 |c(\lambda)|^{-2} d\lambda
$$
where $C_0 > 0$ is a constant, and $c$ is a certain complex function on $\C - \{0\}$.
The required conditions on $A$ are as follows:

\medskip

Making the change of dependent variable $v = A^{1/2} u$, equation (\ref{evaleqn}) becomes
\begin{equation} \label{evaleqn2}
v''(r) = (G(r) - \lambda^2)v(r)
\end{equation}

where the function $G$ is defined by

\begin{equation} \label{gfunction}
G(r) = \frac{1}{4} \left( \frac{A'(r)}{A(r)}\right)^2 + \frac{1}{2}\left(\frac{A'}{A}\right)'(r) - \rho^2
\end{equation}

If the function $G$ tends to $0$ fast enough near infinity, then it is reasonable to expect that equation (\ref{evaleqn2})
above has two linearly independent solutions asymptotic to exponentials $e^{\pm i \lambda r}$ near infinity.
Bloom and Xu show that this is indeed the case \cite{bloomxu1} under the following hypothesis on the function $G$:

\medskip

\noindent (H4) For some $r_0 > 0$, we have
$$
\int_{r_0}^{\infty} r |G(r)| dr < +\infty
$$
and $G$ is bounded on $[r_0, \infty)$.

\medskip

Under hypothesis (H4), for any $\lambda \in \C - \{0\}$, there are unique solutions $\Phi_{\lambda}, \Phi_{-\lambda}$
of equation (\ref{evaleqn}) on $(0, \infty)$ which are asymptotic to exponentials near infinity \cite{bloomxu1},
$$
\Phi_{\pm \lambda}(r) = e^{(\pm i\lambda - \rho)r}(1 + o(1)) \ \hbox{ as } r \to +\infty
$$
The solutions $\Phi_{\lambda}, \Phi_{-\lambda}$ are linearly independent, so, since $\phi_{\lambda} = \phi_{-\lambda}$,
there exists a function $c$ on $\C - \{0\}$ such that
$$
\phi_{\lambda} = c(\lambda) \Phi_{\lambda} + c(-\lambda) \Phi_{-\lambda}
$$
for all $\lambda \in \C - \{0\}$. We will call this function the $c$-function of the hypergroup.
We remark that if the hypergroup $([0, \infty), *)$ is the one arising from
convolution of radial measures on a noncompact rank one symmetric space $G/K$, then this function agrees
with Harish-Chandra's $c$-function only on the half-plane $\{ \Im \lambda \leq 0 \}$ and not on all of $\C$.

\medskip

If we furthermore assume the hypothesis $|\alpha| \neq 1/2$, then Bloom-Xu show that the function $c$ is non-zero for
$\Im \lambda \leq 0, \lambda \neq 0$, and prove the following estimates:

\medskip

There exist constants $C, K > 0$ such that

\begin{align*}
\frac{1}{C} |\lambda| & \leq |c(\lambda)|^{-1} \leq C |\lambda| , \quad \quad  |\lambda|  \leq K \\
\frac{1}{C} |\lambda|^{\alpha + \frac{1}{2}} & \leq |c(\lambda)|^{-1} \leq C |\lambda|^{\alpha + \frac{1}{2}} , \quad |\lambda| \geq K \\
\end{align*}

Moreover they prove the following inversion formula: for any even function $f \in C^{\infty}_c(\R)$,
$$
f(r) = C_0 \int_{0}^{\infty} \hat{f}(\lambda) \phi_{\lambda}(r) |c(\lambda)|^{-2} d\lambda
$$
where $C_0 > 0$ is a constant.

\medskip

It follows that the Plancherel measure $\sigma$ of the hypergroup is supported on $[0, \infty)$, and absolutely
continuous with respect to Lebesgue measure, with density given by $C_0 |c(\lambda)|^{-2}$. Bloom-Xu also
show that the $c$-function is holomorphic on the half-plane $\{ \Im \lambda < 0 \}$.

\medskip

\subsection{The density function of a harmonic manifold}

\medskip

Let $X$ be a simply connected, $n$-dimensional negatively curved harmonic manifold as before, and let $A$ be the density
function of $X$. We check that $A$ is a Chebli-Trimeche function, so that we obtain a commutative hypergroup $([0,\infty), *)$,
and that the conditions of Bloom-Xu are met so that the Plancherel measure is given by $C_0 |c(\lambda)|^{-2} d\lambda$ on
$[0, \infty)$.

\medskip

The function $A(r)$ equals, up to a constant factor, the volume of geodesic spheres $S(x, r)$, which is increasing in $r$ and
tends to infinity as $r$ tends to infinity, so condition (H1) is satisfied. As stated in section 2.2, the function $A'(r)/A(r)$
equals the mean curvature of geodesic spheres $S(x, r)$, which decreases
monotonically to a limit $2\rho$ which is positive (and equals the mean curvature of horospheres), so condition (H2) is satisfied.

\medskip

Fixing a point $x \in X$, for $r > 0$, the density function $A(r)$ is given by the Jacobian of the map $\phi : v \mapsto \exp_x(rv)$ from
the unit tangent sphere $T^1_x X$ to the geodesic sphere $S(x, r)$. Let $T$ be the map $v \mapsto rv$ from the unit tangent sphere
$T^1_x X$ to the tangent sphere of radius $r$, $T^r_x X \subset T_x M$, then $\phi = \exp_x \circ T$, so the Jacobian of $\phi$ is
given by the product of the Jacobians of $T$ and $\exp_x$, hence
$$
A(r) = r^{n-1} B(r)
$$
where the function $B$ is given by
$$
B(r) = \det (D\exp_x)_{rv}
$$
where $v$ is any fixed vector in $T^1_x X$. Since $B$ is independent of the choice of $v$, in particular is the same for vectors $v$ and
$-v$, the function $B$ is even, and $C^{\infty}$ on $\R$ with $B(0) = 1$.
Thus condition (H3) holds for the function $A$, with $\alpha = (n - 2)/2$.

\medskip

The density function $A$ is thus a Chebli-Trimeche function, so we obtain a hypergroup structure on $[0, \infty)$, which we call
the {\it radial hypergroup} of the harmonic manifold $X$ (the reason for this terminology will become clear from the
the following sections).

\medskip

We proceed to check that condition (H4) is satisfied. For this we will need the following theorem of Nikolayevsy:

\medskip

\begin{theorem} \label{exppoly} \cite{nikolayevsky} The density function of a harmonic manifold is an {\it exponential polynomial},
i.e. a function of the form
$$
A(r) = \sum_{i = 1}^k (p_i(r) \cos(\beta_i r) + q_i(r) \sin(\beta_i r)) e^{\alpha_i r}
$$
where $p_i, q_i$ are polynomials and $\alpha_i, \beta_i \in \R$, $i = 1, \dots, k$.
\end{theorem}

\medskip

It will be convenient to rearrange terms and write the density function in the form
\begin{equation} \label{formofa}
A(r) = \sum_{i = 1}^l \sum_{j = 0}^{m_i} f_{ij}(r) r^j e^{\alpha_i r}
\end{equation}
where $\alpha_1 < \alpha_2 < \dots < \alpha_l$, and each $f_{ij}$ is a trigonometric polynomial,
i.e. a finite linear combination of functions of the form $\cos(\beta r)$ and $\sin(\beta r)$, $\beta \in \R$, with
$f_{im_i}$ not identically zero, for $i = 1, \dots, l$.
For an exponential polynomial written in this form, we will call the largest exponent $\alpha_l$ which appears in the
exponentials the {\it exponential degree} of the exponential polynomial.

\medskip

\begin{lemma} \label{densityform} With the density function as above, we have $\alpha_l = 2\rho, m_l = 0$ and $f_{l0} = C$
for some constant $C > 0$. Thus the density function is of the form
$$
A(r) = C e^{2\rho r} + P(r)
$$
where $P$ is an exponential polynomial of exponential degree $\delta < 2\rho$.
\end{lemma}

\medskip

\noindent{\bf Proof:} Since $X$ is CAT(-1), in particular $X$ is Gromov-hyperbolic,
so by a result of Knieper \cite{knieper1} $X$ has {\it purely exponential volume growth}, i.e.
there exists a constant $C > 1$ such that
\begin{equation} \label{pureexp}
\frac{1}{C} \leq \frac{A(r)}{e^{2\rho r}} \leq C
\end{equation}
for all $r \geq 1$. If $\alpha_l < 2\rho$, then $A(r)/e^{2\rho r} \to 0$ as $r \to \infty$, contradicting (\ref{pureexp})
above, so we must have $\alpha_l \geq 2\rho$. On the other hand, if $\alpha_l > 2\rho$, then since
$f_{lm_l}$ is a trigonometric polynomial which is not identically zero, we can choose a sequence $r_m$ tending to infinity
such that $f_{lm_l}(r_m) \to \alpha \neq 0$. Then clearly $A(r_m)/e^{2\rho r_m} \to \infty$, again contradicting (\ref{pureexp}).
Hence $\alpha_l = 2\rho$.

\medskip

Using (\ref{formofa}) and $\alpha_l = 2\rho$, we have
$$
\frac{A'(r)}{A(r)} - 2\rho = \frac{f_{lm_l}'(r) + o(1)}{f_{lm_l}(r) + o(1)}
$$
as $r \to \infty$, thus
\begin{align*}
f_{lm_l}'(r) + o(1) & = (f_{lm_l}(r) + o(1))\left(\frac{A'(r)}{A(r)} - 2\rho\right) \\
                    & \to 0 \\
\end{align*}
as $r \to \infty$ since $f_{lm_l}$ is bounded and $A'(r)/A(r) - 2\rho \to 0$ as $r \to \infty$.
Thus $f_{lm_l}'$ is a trigonometric polynomial which tends to $0$ as $r \to \infty$, so it must be identically
zero, hence $f_{lm_l} = C$ for some non-zero constant $C$.

\medskip

It follows that
$$
A(r) = C r^{m_l} e^{2\rho r} (1 + o(1))
$$
as $r \to \infty$. If $m_l \geq 1$ then $A(r)/e^{2\rho r} \to \infty$ as $r \to \infty$, so we must have $m_l = 0$.
$\diamond$

\medskip

\begin{lemma} \label{h4holds} Condition (H4) holds for the density function $A$, i.e.
$$
\int_{r_0}^{\infty} r |G(r)| dr < +\infty
$$
and $G$ is bounded on $[r_0, \infty)$ for any $r_0 > 0$, where
$$
G(r) = \frac{1}{4} \left( \frac{A'(r)}{A(r)}\right)^2 + \frac{1}{2}\left(\frac{A'}{A}\right)'(r) - \rho^2
$$
\end{lemma}

\medskip

\noindent{\bf Proof:} By the previous lemma, $A(r) = C e^{2\rho r} + P(r)$, where $P$ is an
exponential polynomial of exponential degree $\delta < 2\rho$. We then have
\begin{align*}
\frac{A'(r)}{A(r)} - 2\rho & = \frac{P'(r) - 2\rho P(r)}{C e^{2\rho r} + P(r)} \\
                           & = \frac{Q(r)}{C e^{2\rho r} + P(r)} \\
\end{align*}
where $Q$ is an exponential polynomial of exponential degree less than or equal to $\delta$. Putting $\alpha = (2\rho - \delta)/2$,
it follows that $A'(r)/A(r) - 2\rho = O(e^{-\alpha r})$ as $r \to \infty$.
Differentiating, we obtain
\begin{align*}
\left(\frac{A'}{A}\right)'(r) & = \frac{(C e^{2\rho r} + P(r))Q'(r) - Q(r)(2\rho C e^{2\rho r} + P'(r))}{(C e^{2\rho r} + P(r))^2} \\
                 & = \frac{Q_1(r)}{(C e^{2\rho r} + P(r))^2} \\
\end{align*}
where $Q_1$ is an exponential polynomial of exponential degree less than or equal to $(2\rho + \delta)$. Since the denominator
of the above expression is of the form $k e^{4\rho r} + P_1(r)$ with $P_1$ an exponential polynomial of exponential degree strictly
less than $4\rho$, it follows that $(A'/A)'(r) = O(e^{-\alpha r})$ as $r \to \infty$.

\medskip

Now we can write the function $G$ as
$$
G(r) = \frac{1}{4}\left(\frac{A'(r)}{A(r)} - 2\rho\right)\left(\frac{A'(r)}{A(r)} + 2\rho\right) + \frac{1}{2}\left(\frac{A'}{A}\right)'(r)
$$
Since $(A'(r)/A(r) + 2\rho)$ is bounded, it follows from the previous paragraph that $G(r) = O(e^{-\alpha r})$ as $r \to \infty$.
This immediately implies that condition (H4) holds. $\diamond$

\medskip

In order to apply the result of Bloom-Xu on the Plancherel measure for the hypergroup, it remains to check
that $|\alpha| \neq 1/2$. Since $\alpha = (n - 2)/2$, this means $n \neq 3$. Now the
Lichnerowicz conjecture holds in dimensions $n \leq 5$ (\cite{lichnerowicz1}, \cite{walker1}, \cite{besse1}, \cite{nikolayevsky}), i.e. the only harmonic manifolds in such dimensions
are the rank one symmetric spaces $X = G/K$, for which as mentioned earlier the Jacobi analysis
applies, and the Plancherel measure of the hypergroup is well known to be given by $C_0 |{\bf c}(\lambda)|^{-2} d\lambda$ where
${\bf c}$ is Harish-Chandra's $c$-function. Thus in our case we may as well assume that $X$ has dimension $n \geq 6$, so that
$|\alpha| \neq 1/2$, and we may then apply the results of Bloom-Xu stated in the previous section.

\medskip

\subsection{The spherical Fourier transform}

\medskip

Let $\phi_{\lambda}$ denote as in section 4.1 the unique function on $[0, \infty)$ satisfying $L_R \phi_{\lambda} = -(\lambda^2 + \rho^2) \phi_{\lambda}$
and $\phi_{\lambda}(0) = 1$. For $x \in X$ let $d_x$ denote as before the distance function from the point $x$, $d_x(y) = d(x, y)$.
We define the following eigenfunction of $\Delta$ radial around $x$:
$$
\phi_{\lambda, x} := \phi_{\lambda} \circ d_x
$$
The uniqueness of $\phi_{\lambda}$ as an eigenfunction of $L_R$ with eigenvalue $-(\lambda^2 + \rho^2)$ and taking the value $1$ at $r = 0$
immediately implies the following lemma:

\medskip

\begin{lemma} \label{uniqueeigenfn} The function $\phi_{\lambda, x}$ is the unique eigenfunction $f$ of $\Delta$ on $X$ with
eigenvalue $-(\lambda^2 + \rho^2)$ which is radial around $x$ and satisfies $f(x) = 1$.
\end{lemma}

\medskip

Note that for $\lambda \in \R$, the functions $\phi_{\lambda, x}$ are bounded. Let $dvol$ denote the Riemannian volume measure on $X$.

\medskip

\begin{definition} Let $f \in L^1(X, dvol)$ be radial around the point $x \in X$. We define the spherical Fourier transform
of $f$ by
$$
\hat{f}(\lambda) := \int_X f(y) \phi_{\lambda, x}(y) dvol(y)
$$
for $\lambda \in \R$.
\end{definition}

\medskip

For $f$ a function on $X$ radial around the point $x$,
let $f = u \circ d_x$ where $u$ is a function on $[0, \infty)$, then evaluating the integral over $X$ in
geodesic polar coordinates gives
$$
\int_X |f(y)| dvol(y) = \int_{0}^{\infty} |u(r)| A(r) dr
$$
thus $f \in L^1(X)$ if and only if $u \in L^1([0, \infty), A(r) dr)$. In that case, again integrating in polar
coordinates gives
$$
\hat{f}(\lambda) = \int_{0}^{\infty} u(r) \phi_{\lambda}(r) A(r) dr = \hat{u}(\lambda)
$$
where $\hat{u}$ is the hypergroup Fourier transform of the function $u$. Moreover $f \in C^{\infty}_c(X)$ if and only if
$u$ extends to an even function on $\R$ such that $u \in C^{\infty}_c(\R)$. Applying the Fourier inversion formula of Bloom-Xu for the radial
hypergroup stated in section 4.1 to the function $u$ then leads immediately to the following inversion formula for radial functions:

\medskip

\begin{theorem} \label{radialinversion} Let $f \in C^{\infty}_c(X)$ be radial around the point $x \in X$. Then
$$
f(y) = C_0 \int_{0}^{\infty} \hat{f}(\lambda) \phi_{\lambda, x}(y) |c(\lambda)|^{-2} d\lambda
$$
for all $y \in X$. Here $c$ denotes the $c$-function of the radial hypergroup and $C_0 > 0$ is a constant.
\end{theorem}

\medskip

\noindent{\bf Proof:} As shown in the previous section, all the hypotheses required to apply the inversion formula of Bloom-Xu are satisfied,
hence
$$
u(r) = C_0 \int_{0}^{\infty} \hat{u}(\lambda) \phi_{\lambda}(r) |c(\lambda)|^{-2} d\lambda
$$
Since $f = u \circ d_x$, this gives
\begin{align*}
f(y) & = u(d_x(y)) \\
     & = C_0 \int_{0}^{\infty} \hat{u}(\lambda) \phi_{\lambda}(d_x(y)) |c(\lambda)|^{-2} d\lambda \\
     & = C_0 \int_{0}^{\infty} \hat{f}(\lambda) \phi_{\lambda,x}(y) |c(\lambda)|^{-2} d\lambda \\
\end{align*}
$\diamond$

\medskip

The Plancherel theorem for the radial hypergroup leads to the following:

\medskip

\begin{theorem} Let $L^2_x(X, dvol)$ denote the closed subspace of $L^2(X)$ consisting of those functions
in $L^2(X)$ which are radial around the point $x$. For $f \in L^1(X, dvol) \cap L^2_x(X, dvol)$, we have
$$
\int_X |f(y)|^2 dvol(y) = C_0 \int_{0}^{\infty} |\hat{f}(\lambda)|^2 |c(\lambda)|^{-2} d\lambda
$$
The spherical Fourier transform $f \mapsto \hat{f}$ extends to an isometry from $L^2_x(X, dvol)$ onto
$L^2([0, \infty), C_0 |c(\lambda)|^{-2} d\lambda)$.
\end{theorem}

\medskip

\noindent{\bf Proof:} The map $u \mapsto f = u \circ d_x$ defines
an isometry of $L^2([0, \infty), A(r) dr)$ onto $L^2(X, dvol)_x$, which maps $L^1([0, \infty), A(r) dr) \cap L^2([0, \infty), A(r) dr)$
onto $L^1(X, dvol) \cap L^2_x(X, dvol)$. The statements of the theorem then follow from the Levitan-Plancherel theorem
for the radial hypergroup and from the fact that the Plancherel measure is supported on $[0, \infty)$,
given by $C_0 |c(\lambda)|^{-2} d\lambda$. $\diamond$

\medskip

\section{Fourier inversion and Plancherel theorem}

\medskip

We proceed to the analysis of non-radial functions on $X$. Our definition of Fourier transform will
depend on the choice of a basepoint $x \in X$.

\medskip

\begin{definition} Let $x \in X$. For $f \in C^{\infty}_c(X)$, the Fourier transform of $f$ based at the
point $x$ is the function
on $\C \times \partial X$ defined by
$$
\tilde{f}^x(\lambda, \xi) = \int_X f(y) e^{(-i\lambda - \rho)B_{\xi,x}(y)} dvol(y)
$$
for $\lambda \in \C, \xi \in \partial X$. Here as before $B_{\xi,x}$ denotes the Busemann function at
$\xi$ based at $x$ such that $B_{\xi,x}(x) = 0$.
\end{definition}

\medskip

Using the formula
$$
B_{\xi, x} = B_{\xi,o} - B_{\xi, o}(x)
$$
for points $o, x \in X$, we obtain the following relation between the Fourier transforms
based at two different basepoints $o, x \in X$:

\medskip

\begin{equation} \label{fourierbasept}
\tilde{f}^x(\lambda, \xi) = e^{(i\lambda + \rho)B_{\xi, o}(x)} \tilde{f}^o(\lambda, \xi)
\end{equation}

\medskip

The key to passing from the inversion formula for radial functions of section 4.3 to an inversion formula for non-radial
functions will be a formula expressing the radial eigenfunctions $\phi_{\lambda,x}$ as an integral with respect to
$\xi \in \partial X$ of the eigenfunctions $e^{(i\lambda - \rho)B_{\xi,x}}$ (Theorem \ref{phipoisson}). This will be the analogue of the well-known
formulae for rank one symmetric spaces $G/K$ and harmonic $NA$ groups expressing the radial eigenfunctions $\phi_{\lambda, x}$ as matrix
coefficients of representations of $G$ on $L^2(K/M)$ and $NA$ on $L^2(N)$ respectively. We first need to define the {\it visibility measures}
on the boundary $\partial X$:

\medskip

Given a point $x \in X$, let $\lambda_x$ be normalized Lebesgue measure on the unit tangent sphere $T^1_x X$,
i.e. the unique probability measure on $T^1_x X$ invariant under the orthogonal group of the tangent space
$T_x M$. For $v \in T^1_x X$, let $\gamma_v : [0, \infty) \to X$ be the unique geodesic ray with initial velocity $v$.
Then we have a homeomorphism $p_x : T^1_x X \to \partial X, v \mapsto \gamma_v(\infty)$. The
visibility measure on $\partial X$ (with respect to the basepoint $x$)
is defined to be the push-forward $(p_x)_* \lambda_x$ of $\lambda_x$ under the map $p_x$; for
notational convenience, we will however denote the visibility measure on $\partial X$ by the same symbol $\lambda_x$.

\medskip

In \cite{knieperpeyerimhoff1}, it is shown that the visibility measures $\lambda_x, x \in X$ are mutually absolutely
continuous and their Radon-Nikodym derivatives are given by
$$
\frac{d\lambda_y}{d\lambda_x}(\xi) = e^{-2\rho B_{\xi, x}(y)}
$$
The above formula for Radon-Nikodym derivatives leads in a standard way to the following estimate for visual measures
of visual balls (in the literature the estimate is usually proved for "shadows" of balls in $X$ and called
Sullivan's Shadow Lemma):

\medskip

\begin{lemma} \label{shadow} There is a constant $C > 1$ such that
$$
\frac{1}{C} \epsilon^{2\rho} \leq \lambda_x(B_x(\xi, \epsilon)) \leq C \epsilon^{2\rho}
$$
for all $x \in X, \xi \in \partial X, \epsilon > 0$.
\end{lemma}

\medskip

\noindent{\bf Proof:} Given $x \in X, \xi \in \partial X$ and $\epsilon > 0$, choose $y$ on the geodesic ray
$[x, \xi)$ at distance $r = \log(1/\epsilon)$ from $x$. Using the formula for the Radon-Nikodym
derivative of $\lambda_x$ with respect to $\lambda_y$, we have
$$
\lambda_x(B_x(\xi, \epsilon)) = \int_{B_x(\xi, \epsilon)} e^{-2\rho B(x, y, \eta)} d\lambda_y(\eta)
$$
On the other hand by Lemma \ref{changebasept} we know that $|B(x, y, \eta) - r| \leq 2$ for any $\eta \in B_x(\xi, \epsilon)$.
Using $e^{-2\rho r} = \epsilon^{2\rho}$, it follows that
$$
\frac{1}{C} \lambda_y(B_x(\xi, \epsilon)) \epsilon^{2\rho} \leq \lambda_x(B_x(\xi, \epsilon)) \leq C \lambda_y(B_x(\xi, \epsilon)) \epsilon^{2\rho}
$$
for some constant $C$. Since $\lambda_y$ is a probability measure,
this gives the required upper bound
$$
\lambda_x(B_x(\xi, \epsilon)) \leq C \epsilon^{2\rho}
$$
for some constant $C > 1$.

\medskip

For the lower bound, it suffices to give a lower bound on $\lambda_y(B_x(\xi, \epsilon))$ by a positive constant. By Lemma \ref{changebasept},
there is a universal constant $\delta_0 > 0$ such that $B_y(\xi, \delta_0) \subset B_x(\xi, \epsilon)$, so
$$
\lambda_y(B_x(\xi, \epsilon)) \geq \lambda_y(B_y(\xi, \delta_0)).
$$
Now the sectional curvature of $X$ is bounded below, $-b^2 \leq K$ (\cite{knieper2}, Cor. 2.12 in section 1.2).
Identifying $\partial X$ with $T^1_y X$ via the map $p_y : T^1_y X \to \partial X$ and then identifying
$T^1_y X$ with the standard sphere $S^{n-1}$ via an isometry $T : \R^n \to T^1_y X$, the measure $\lambda_y$
corresponds to Lebesgue measure on $S^{n-1}$, and it follows from Lemma \ref{anglecomp}
that the image of the ball $B_y(\xi, \delta_0)$ in $S^{n-1}$ contains a ball in $S^{n-1}$ of radius $\delta^b_0$
with respect to the chordal metric, which has Lebesgue measure $\alpha > 0$ where $\alpha$ only depends on $\delta_0$, hence
$$
\lambda_y(B_y(\xi, \delta_0)) \geq \alpha
$$
and the Lemma follows.
$\diamond$

\medskip

\begin{lemma} \label{intconv} Let $x \in X$ and $\xi \in \partial X$. Then the integral
$$
c(\lambda, x, \xi) := \int_{\partial X} \rho_x(\xi, \eta)^{2(i\lambda - \rho)} d\lambda_x(\eta)
$$
converges for $\Im \lambda < 0$ and is holomorphic in $\lambda$ on the half-plane $\{ \Im \lambda < 0 \}$.
\end{lemma}

\medskip

\noindent{\bf Proof:} Let $\lambda = \sigma - i\tau$, where $\sigma \in \R$ and $\tau > 0$. Then
\begin{align*}
\int_{\partial X} |\rho_x(\xi, \eta)^{2(i\lambda - \rho)}| d\lambda_x(\eta) & = \int_{\partial X} \rho_x(\xi, \eta)^{2(\tau - \rho)} d\lambda_x(\eta) \\
                                                                            & = \int_{0}^{\infty} \lambda_x(\{ \eta \ | \ \rho_x(\xi, \eta)^{2(\tau - \rho)} > t \}) dt \\
\end{align*}

If $\tau \geq \rho$ then the set $\{ \eta \ | \ \rho_x(\xi, \eta)^{2(\tau - \rho)} > t \}$ is empty for $t > 1$ (since the visual metrics have diameter one),
and so the last integral reduces to an integral over $[0,1]$,
which is bounded above by one
since $\lambda_x$ is a probability measure.

\medskip

For $0 < \tau < \rho$ using Lemma \ref{shadow} and the fact that $\lambda_x$ is a probability measure we have
\begin{align*}
\int_{0}^{\infty} \lambda_x(\{ \eta \ | \ \rho_x(\xi, \eta)^{2(\tau - \rho)} > t \}) dt & \leq 1 + \int_{1}^{\infty} \lambda_x(B_x(\xi, (1/t)^{1/2(\rho - \tau)})) dt \\
& \leq 1 + C \int_{1}^{\infty} \left(\frac{1}{t}\right)^{\frac{2\rho}{2(\rho - \tau)}} dt \\
& < +\infty \\
\end{align*}

Thus the integral defining $c(\lambda, x, \xi)$ converges for $\Im \lambda < 0$. That $c(\lambda, x, \xi)$ is holomorphic
in $\lambda$ follows from Morera's theorem. $\diamond$

\medskip

\begin{lemma} \label{radialisation} Let $x \in X$ and $\xi \in \partial X$. Then for all $\lambda \in \C$,
$$
\phi_{\lambda, x} = M_x (e^{(i\lambda - \rho)B_{\xi,x}})
$$
(where $M_x$ is the radialisation operator around the point $x$). In particular, $\phi_{\lambda, x}(y)$ is entire in
$\lambda$ for fixed $y \in X$, and is real and positive for $\lambda$ such that $(i\lambda - \rho)$ is real and positive.
\end{lemma}

\medskip

\noindent{\bf Proof:} Since the function $e^{(i\lambda - \rho)B_{\xi,x}}$ is an eigenfunction of the
Laplacian $\Delta$ with eigenvalue $-(\lambda^2 + \rho^2)$ and the operator $M_x$ commutes
with $\Delta$, the function $f = M_x (e^{(i\lambda - \rho)B_{\xi,x}})$ is also an eigenfunction of
$\Delta$ for the eigenvalue $-(\lambda^2 + \rho^2)$. Since $f$ is radial around $x$ and $f(x) = 1$, it follows
from Lemma \ref{uniqueeigenfn} that $f = \phi_{\lambda, x}$. $\diamond$

\medskip

\begin{prop} \label{coincide} Let $f \in C^{\infty}_c(X)$ be radial around the point $x \in X$. Then the Fourier transform of $f$
based at $x$ coincides with the spherical Fourier transform,
$$
\tilde{f}^x(\lambda, \xi) = \hat{f}(\lambda)
$$
for all $\lambda \in \C, \xi \in \partial X$.
\end{prop}

\medskip

\noindent{\bf Proof:} Let $f = u \circ d_x$ where $u \in C^{\infty}_c(\R)$. By Lemma \ref{radialisation} above,
$$
\phi_{\lambda}(r) = \phi_{-\lambda}(r) = \int_{S(x,r)} e^{(-i\lambda - \rho)B_{\xi,x}(y)} d\sigma^r(y)
$$
where $\sigma^r$ is normalized surface area measure on the geodesic sphere $S(x,r)$.
Evaluating the integral defining $\tilde{f}^x$ in geodesic polar coordinates centered at $x$ we have
\begin{align*}
\tilde{f}^x(\lambda, \xi) & = \int_{0}^{\infty} \int_{S(x,r)} f(y) e^{(-i\lambda - \rho)B_{\xi,x}(y)} d\sigma^r(y) A(r) dr \\
                          & = \int_{0}^{\infty} u(r) \phi_{\lambda}(r) A(r) dr \\
                          & = \hat{f}(\lambda) \\
\end{align*}
$\diamond$

\medskip

\begin{prop} \label{cfunction} Let $c$ be the $c$-function of the radial hypergroup of $X$ and let $\Im \lambda < 0$. Then:

\medskip

\noindent (1) We have
$$
\lim_{r \to \infty} \frac{\phi_{\lambda}(r)}{e^{(i\lambda - \rho)r}} = c(\lambda)
$$

\medskip

\noindent (2) We have
$$
c(\lambda) = \int_{\partial X} \rho_x(\xi, \eta)^{2(i\lambda - \rho)} d\lambda_x(\eta)
$$
for any $x \in X, \xi \in \partial X$. In particular the integral $c(\lambda, x, \xi)$ above is independent of the choice of $x, \xi$.
\end{prop}

\medskip

\noindent{\bf Proof:} (1): For $\Im \lambda < 0$, using $\phi_{\lambda} = c(\lambda) \Phi_{\lambda} + c(-\lambda) \Phi_{-\lambda}$ and
$\Phi_{\pm \lambda}(r) = e^{(\pm i \lambda - \rho)r}(1 + o(1))$ as $r \to \infty$, we have
\begin{align*}
\frac{\phi_{\lambda}(r)}{e^{(i\lambda - \rho)r}} & = c(\lambda)(1 + o(1)) + c(-\lambda)e^{-2i\lambda r}(1 + o(1)) \\
                                                 & \to c(\lambda) \\
\end{align*}
as $r \to \infty$.

\medskip

\noindent (2): Let $\lambda = it$ where $t \leq -\rho$, so that $\mu := i\lambda - \rho \geq 0$. Fix $x \in X$ and $\xi \in \partial X$.
For $\eta \in \partial X$, let $y(\eta, r) \in X$ denote the point
on the geodesic $[x, \eta)$ at a distance $r$ from $x$. Then the normalized surface area measure on the geodesic sphere $S(x, r)$ is given by the
push-forward of $\lambda_x$ under the map $\eta \mapsto y(\eta, r)$, so by Lemma \ref{radialisation}
$$
\frac{\phi_{\lambda}(r)}{e^{(i\lambda - \rho)r}} = \int_{\partial X} e^{(i\lambda - \rho)(B_{\xi, x}(y(\eta, r)) - r)} d\lambda_x(\xi)
$$

We will apply the dominated convergence theorem to evaluate the limit of the above integral as $r \to \infty$. First note that
by Lemma \ref{busemasymptotic}, for any $\eta$ not equal to $\xi$,
$$
B_{\xi, x}(y(\eta, r)) - r \to -2(\xi|\eta)_x
$$
as $r \to \infty$, so the integrand converges a.e. as $r \to \infty$,
$$
e^{(i\lambda - \rho)(B_{\xi, x}(y(\eta, r)) - r)} \to \rho_x(\xi, \eta)^{2(i\lambda - \rho)}
$$
(using $\rho_x(\xi, \eta) = e^{-(\xi|\eta)_x}$).

\medskip

Now by Lemma \ref{buseform}, using $\mu \geq 0$ we can estimate the integrand:
\begin{align*}
e^{\mu (B_{\xi, x}(y(\eta, r)) - r)} & = (\cosh r - \sinh r \cos(\theta_x(\xi, y(x,r))))^{\mu} e^{-\mu r} \\
                                     & = (e^{-r} + 2 \sinh r \sin^2(\theta_x(\xi, y(x, r))/2))^{\mu} e^{-\mu r} \\
                                     & \leq (e^{-r} + (e^r - e^{-r}) \cdot 1)^{\mu} e^{-\mu r} \\
                                     & = 1 \\
\end{align*}

So dominated convergence applies and we conclude that
$$
\frac{\phi_{\lambda}(r)}{e^{(i\lambda - \rho)r}} \to \int_{\partial X} \rho_x(\xi, \eta)^{2(i\lambda - \rho)} d\lambda_x(\eta)
$$
as $r \to \infty$ for $\lambda = it, t \leq -\rho$.

\medskip

It follows from part (1) of the proposition that $c(\lambda)$ equals the integral above, which is $c(\lambda, x, \xi)$,
for $\lambda = it, t \leq -\rho$.
Since $c(\lambda)$ and $c(\lambda, x, \xi)$ are holomorphic for $\Im \lambda < 0$, they must then be equal for all $\lambda$
with $\Im \lambda < 0$. $\diamond$

\medskip

\noindent{\bf Remark.} Part (2) of the above proposition is the analogue of the well-known
integral formula for Harish-Chandra's {\bf c}-function (formula (18) in \cite{helgason1}, pg. 108).

\medskip

For $\lambda \in \C$ and $x \in X$, define the function $\tilde{\phi}_{\lambda, x}$ on $X$ by
$$
\tilde{\phi}_{\lambda,x}(y) = \int_{\partial X} e^{(i\lambda - \rho)B_{\xi, x}(y)} d\lambda_x(y)
$$
It follows from the above equation that $\tilde{\phi}_{\lambda, x}(y)$ is entire in $\lambda$ for fixed $y \in X$,
and is real and positive for $\lambda$ such that $(i\lambda - \rho)$ is real and positive. Moreover, by
Proposition \ref{helgasonkernel}, the function $\tilde{\phi}_{\lambda, x}$ is an eigenfunction of the Laplacian $\Delta$ with
eigenvalue $-(\lambda^2 + \rho^2)$, and $\tilde{\phi}_{\lambda, x}(x) = 1$.

\medskip

\begin{lemma} \label{almostradial} Let $\lambda = it, t < -\rho$, and let $x \in X$. Given $\epsilon > 0$,
there exists $R = R(\epsilon)$ such that for any $y \in X$ with $r = d(x,y) \geq R$, we have
$$
\left|\frac{\tilde{\phi}_{\lambda,x}(y)}{e^{(i\lambda - \rho)r}} - c(\lambda)\right| \leq \epsilon
$$
\end{lemma}

\medskip

\noindent{\bf Proof:} Let $\mu = i\lambda - \rho > 0$. We may assume $\epsilon > 0$ is small enough so that $|e^{\mu \delta} - 1| \leq 2\mu \epsilon$
for $0 \leq \delta \leq \epsilon$.

\medskip

Let $y \in X$ with $r = d(x,y) > 0$, then the geodesic segment $[x, y]$ extends uniquely to a
geodesic ray $[x, y^+)$ for some $y^+ \in \partial X$. By Lemma \ref{busemasymptotic}, there exists $R = R(\epsilon)$ such that
if $r \geq R$ and $\xi \in \partial X$ is such that $\rho_x(\xi, y^+) \geq \epsilon$ then
$$
|(B_{\xi,x}(y) - r) + 2(\xi|y^+)_x| < \epsilon
$$
and hence we can write
$$
\frac{\exp^{\mu B_{\xi, x}(y)}}{e^{\mu r}} = \rho_x(\xi, y^+)^{2\mu} e^{\mu \delta(\xi)}
$$
where $\delta = \delta(\xi)$ satisfies $|\delta(\xi)| \leq \epsilon$ for $\rho_x(\xi, y^+) \geq \epsilon$.
Now
\begin{align*}
\frac{\tilde{\phi}_{\lambda, x}(y)}{e^{\mu r}} & = \int_{\partial X} \frac{e^{\mu B_{\xi, x}(y)}}{e^{\mu r}} d\lambda_x(\xi) \\
           & = \int_{\partial X - B_x(y^+, \epsilon)} \frac{e^{\mu B_{\xi, x}(y)}}{e^{\mu r}} d\lambda_x(\xi) +
               \int_{B_x(y^+, \epsilon)} \frac{e^{\mu B_{\xi, x}(y)}}{e^{\mu r}} d\lambda_x(\xi) \\
           & = I_1 + I_2 \\
\end{align*}
say. We estimate the two integrals $I_1, I_2$ appearing above separately.

\medskip

For the second integral $I_2$, since $\mu \geq 0$, the same estimate as in the proof of Proposition \ref{cfunction}
gives $e^{\mu B_{\xi, x}(y)}/e^{\mu r} \leq 1$, hence by Lemma \ref{shadow} we obtain
$$
|I_2| \leq C \epsilon^{2\rho}
$$
for some constant $C > 0$.

\medskip

For the first integral $I_1$, by Lemma \ref{cfunction} since $c(\lambda) = c(\lambda, x, y^+)$ we have
\begin{align*}
|I_1 - c(\lambda)| & = \left|\int_{\partial X - B_x(y^+, \epsilon)} \frac{e^{\mu B_{\xi, x}(y)}}{e^{\mu r}} d\lambda_x(\xi)
 - \int_{\partial X} \rho_x(\xi, y^+)^{2\mu} d\lambda_x(\xi)\right| \\
                   & \leq \int_{\partial X - B_x(y^+, \epsilon)} \rho_x(\xi, y^+)^{2\mu} |e^{\mu \delta(\xi)} - 1| d\lambda_x(\xi)
                   + \int_{B_x(y^+, \epsilon)} \rho_x(\xi, y^+)^{2\mu} d\lambda_x(\xi) \\
                   & \leq 2\mu \epsilon + C \epsilon^{2\rho} \\
\end{align*}
(where we have used $\rho_x(\xi, y^+) \leq 1$, $\mu > 0$, and Lemma \ref{shadow}).

\medskip

Putting together the estimates for $I_1$ and $I_2$ gives
$$
\left|\frac{\tilde{\phi}_{\lambda,x}(y)}{e^{(i\lambda - \rho)r}} - c(\lambda)\right| \leq 2\mu \epsilon + 2C \epsilon^{2\rho}
$$
for all $y \in X$ such that $r = d(x,y) \geq R(\epsilon)$. Since $\epsilon > 0$ was arbitrary, the Lemma follows.
$\diamond$

\medskip

\begin{lemma} \label{ratio} Let $\lambda = it, t < -\rho$ and let $x \in X$. Given $\epsilon > 0$,
there exists $R = R(\epsilon)$ such that for any $y \in X$ with $r = d(x,y) \geq R$, we have
$$
1 - \epsilon \leq \frac{\tilde{\phi}_{\lambda, x}(y)}{\phi_{\lambda, x}(y)} \leq 1 + \epsilon
$$
(note that for the given value of $\lambda$ both functions appearing above are real and positive).
\end{lemma}

\medskip

\noindent{\bf Proof:} Note $c(\lambda) > 0$ by Lemma \ref{cfunction}, so by Lemmas \ref{cfunction} and \ref{almostradial}
there exists $R = R(\epsilon)$ such that for $y \in X$ with $r = d(x,y) \geq R$ we have
$$
(1 - \epsilon) c(\lambda) \leq \frac{\tilde{\phi}_{\lambda,x}(y)}{e^{(i\lambda - \rho)r}} \leq (1 + \epsilon) c(\lambda)
$$
and
$$
(1 - \epsilon) c(\lambda) \leq \frac{\phi_{\lambda,x}(y)}{e^{(i\lambda - \rho)r}} \leq (1 + \epsilon) c(\lambda)
$$
The above two inequalities imply
$$
\frac{1 - \epsilon}{1 + \epsilon} \leq \frac{\tilde{\phi}_{\lambda, x}(y)}{\phi_{\lambda, x}(y)} \leq \frac{1 + \epsilon}{1 - \epsilon}
$$
Since $\epsilon > 0$ was arbitrary, the Lemma follows. $\diamond$

\medskip

We can finally prove the required formula expressing $\phi_{\lambda, x}$ as an integral over $\xi \in \partial X$ of the
functions $e^{(i\lambda - \rho)B_{\xi,x}}$:

\medskip

\begin{theorem} \label{phipoisson} Let $\lambda \in \C$ and $x \in X$. Then
$$
\phi_{\lambda, x}(y) = \int_{\partial X} e^{(i\lambda - \rho)B_{\xi,x}(y)} d\lambda_x(\xi)
$$
for all $y \in X$ (i.e. $\phi_{\lambda, x} = \tilde{\phi}_{\lambda, x}$).
\end{theorem}

\medskip

\noindent{\bf Proof:} We first assume $\lambda$ is of the form $\lambda = it, t < -\rho$.

\medskip

For notational convenience, write $\phi = \phi_{\lambda,x}$ and $\tilde{\phi} = \tilde{\phi}_{\lambda,x}$, then $\phi, \tilde{\phi}$
are positive functions. Let $u$ be the function $u = \tilde{\phi}/\phi$.

\medskip

Using the fact that $\phi, \tilde{\phi}$ are both
eigenfunctions of $\Delta$ for the eigenvalue $-(\lambda^2 + \rho^2)$ and the classical formula
$$
\Delta(\phi u) = (\Delta \phi) u + 2 < \nabla \phi, \nabla u > + \phi \Delta u
$$
we have
\begin{align*}
-(\lambda^2 + \rho^2) \tilde{\phi} & = \Delta(\phi u) \\
                                   & = -(\lambda^2 + \rho^2) \phi u + 2 < \nabla \phi, \nabla u > +  \phi \Delta u \\
                                   & = -(\lambda^2 + \rho^2) \tilde{\phi} + 2 < \nabla \phi, \nabla u > +  \phi \Delta u \\
\end{align*}

Dividing the above equation by $\phi$, we obtain
$$
\Delta' u = 0
$$
where $\Delta'$ is the second-order differential operator defined by
$$
\Delta' f = \Delta f + 2 < \frac{\nabla \phi}{\phi}, \nabla f >
$$
for $f \in C^{\infty}(X)$.

\medskip

Now given $\epsilon > 0$, let $R = R(\epsilon) > 0$ be as given by Lemma \ref{ratio}. Then for any $r \geq R$, by Lemma \ref{ratio} we have
$$
(1 - \epsilon) \leq u(y) \leq (1 + \epsilon)
$$
for all $y$ on the geodesic sphere $S(x, r)$. Since the operator $\Delta'$ is an elliptic second order differential operator without
any zeroth order term, and $\Delta' u = 0$, the maximum principle applies to $u$ and $-u$ (\cite{evans}, Theorem 1 of section 6.4.1),
so we obtain
$$
(1 - \epsilon) \leq u(y) \leq (1 + \epsilon)
$$
for all $y$ in the geodesic ball $B(x, r)$. Since this is true for all $r \geq R$, the above inequality holds in fact for all $y \in X$.
Since $\epsilon > 0$ is arbitrary, it follows that $u = 1$ on $X$. Thus
$$
\tilde{\phi}_{\lambda,x}(y) = \phi_{\lambda, x}(y)
$$
for all $y \in X$, and $\lambda$ of the form $\lambda = it, t < -\rho$. Since both sides of the above equation are entire in $\lambda$
for fixed $y$, it follows that $\tilde{\phi}_{\lambda,x} = \phi_{\lambda, x}$ on $X$ for all $\lambda$ as required. $\diamond$

\medskip

We can now prove the Fourier inversion formula:

\medskip

\begin{theorem} \label{fourierinversion} Fix a basepoint $o \in X$. Then for $f \in C^{\infty}_c(X)$ we have
$$
f(x) = C_0 \int_{0}^{\infty} \int_{\partial X} \tilde{f}^o(\lambda, \xi) e^{(i\lambda - \rho)B_{\xi, o}(x)} d\lambda_o(\xi) |c(\lambda)|^{-2} d\lambda
$$
for all $x \in X$ (where $C_0 > 0$ is a constant).
\end{theorem}

\medskip

\noindent{\bf Proof:} Given $f \in C^{\infty}_c(X)$ and $x \in X$, the function $M_x f$ is in $C^{\infty}_c(X)$, is radial
around the point $x$ and satisfies $(M_x f)(x) = f(x)$. By Theorem \ref{radialinversion} applied to the function $M_x f$ we have
\begin{align*} \label{fx}
f(x) = (M_x f)(x) & = C_0 \int_{0}^{\infty} \widehat{M_x f}(\lambda) \phi_{\lambda, x}(x) |c(\lambda)|^{-2} d\lambda \\
                  & = C_0 \int_{0}^{\infty} \widehat{M_x f}(\lambda) |c(\lambda)|^{-2} d\lambda \\
\end{align*}
(since $\phi_{\lambda, x}(x) = 1$). Now using the formal self-adjointness of the operator $M_x$, Theorem \ref{phipoisson},
the fact that $\phi_{\lambda,x}$ is radial around $x$ and $\phi_{\lambda, x} = \phi_{-\lambda,x}$ we obtain
\begin{align*}
\widehat{M_x f}(\lambda) & = \int_X (M_x f)(y) \phi_{-\lambda, x}(y) dvol(y) \\
                     & = \int_X f(y) (M_x \phi_{-\lambda, x})(y) dvol(y) \\
                     & = \int_X f(y) \phi_{-\lambda, x}(y) dvol(y) \\
                     & = \int_X f(y) \left( \int_{\partial X} e^{(-i\lambda - \rho)B_{\xi, x}(y)} d\lambda_x(\xi) \right) dvol(y) \\
                     & = \int_{\partial X} \left(\int_X f(y) e^{(-i\lambda - \rho)B_{\xi, x}(y)} dvol(y)\right) d\lambda_x(\xi) \\
                     & = \int_{\partial X} \tilde{f}^x(\lambda, \xi) d\lambda_x(\xi) \\
\end{align*}

Using the relations
$$
\tilde{f}^x(\lambda, \xi) = e^{(i\lambda + \rho)B_{\xi, o}(x)} \tilde{f}^o(\lambda, \xi)
$$
and
$$
\frac{d\lambda_x}{d\lambda_o}(\xi) = e^{-2\rho B_{\xi, o}(x)}
$$
we get
\begin{align*}
\widehat{M_x f}(\lambda) & = \int_{\partial X} e^{(i\lambda + \rho)B_{\xi, o}(x)} \tilde{f}^o(\lambda, \xi) e^{-2\rho B_{\xi, o}(x)} d\lambda_o(\xi) \\
                     & = \int_{\partial X} \tilde{f}^o(\lambda, \xi) e^{(i\lambda - \rho)B_{\xi, o}(x)} d\lambda_o(\xi) \\
\end{align*}

Substituting this last expression for $\widehat{M_x f}(\lambda)$ in the equation $f(x) = C_0 \int_{0}^{\infty} \widehat{M_x f}(\lambda) |c(\lambda)|^{-2} d\lambda$
gives
$$
f(x) = C_0 \int_{0}^{\infty} \int_{\partial X} \tilde{f}^o(\lambda, \xi) e^{(i\lambda - \rho)B_{\xi, o}(x)} d\lambda_o(\xi) |c(\lambda)|^{-2} d\lambda
$$
as required.
$\diamond$

\medskip

The Fourier inversion formula leads immediately to a Plancherel theorem:

\medskip

\begin{theorem} \label{plancherel} Fix a basepoint $o \in X$. For $f, g \in C^{\infty}_c(X)$, we have
$$
\int_X f(x) \overline{g(x)} dvol(x) = C_0 \int_{0}^{\infty} \int_{\partial X} \tilde{f}^o(\lambda, \xi) \overline{\tilde{g}^o(\lambda, \xi)} d\lambda_o(\xi) |c(\lambda)|^{-2} d\lambda
$$
where $C_0$ is the constant appearing in the Fourier inversion formula.
The Fourier transform $f \mapsto \tilde{f}^o$ extends to an isometry of $L^2(X, dvol)$ into
$L^2([0, \infty) \times \partial X, C_0 |c(\lambda)|^{-2} d\lambda d\lambda_o(\xi))$.
\end{theorem}

\medskip

\noindent{\bf Proof:} Applying the Fourier inversion formula to the function $g$ gives
\begin{align*}
\int_X f(x) \overline{g(x)} dvol(x) & = C_0 \int_X f(x) \left(\int_{0}^{\infty} \int_{\partial X} \overline{\tilde{g}^o(\lambda, \xi)} e^{(-i\lambda - \rho)B_{\xi, o}(x)} d\lambda_o(\xi) |c(\lambda)|^{-2} d\lambda\right) dvol(x) \\
    & = C_0 \int_{0}^{\infty} \int_{\partial X} \left(\int_X f(x) e^{(-i\lambda - \rho)B_{\xi, o}(x)} dvol(x)\right) \overline{\tilde{g}^o(\lambda, \xi)} d\lambda_o(\xi) |c(\lambda)|^{-2} d\lambda \\
    & = C_0 \int_{0}^{\infty} \int_{\partial X} \tilde{f}^o(\lambda, \xi) \overline{\tilde{g}^o(\lambda, \xi)} d\lambda_o(\xi) |c(\lambda)|^{-2} d\lambda \\
\end{align*}
Taking $f = g$ gives that the Fourier transform preserves $L^2$ norms,
$$
||f||_2 = ||\tilde{f}^o||_2
$$
for all $f \in C^{\infty}_c(X)$, from which it follows from a standard argument that
the Fourier transform extends to an isometry of $L^2(X, dvol)$ into
$L^2([0, \infty) \times \partial X, C_0 |c(\lambda)|^{-2} d\lambda d\lambda_o(\xi))$.
$\diamond$

\medskip

\section{The convolution algebra of radial functions}

\medskip

Fix a basepoint $o \in X$. We define a notion of convolution with radial functions as follows:

\medskip

For a function $f$ radial around the point $o$, let $f = u \circ d_o$, where $u$ is a function on $\R$. For $x \in X$,
the {\it x-translate} of $f$ is defined to be the function
$$
\tau_x f = u \circ d_x
$$
Note that if $f \in L^1(X, dvol)$, then evaluating integrals in geodesic polar coordinates centered at $o$ and $x$ gives
$$
||f||_1 = \int_{0}^{\infty} |u(r)| A(r) dr = ||\tau_x f||_1
$$

\begin{definition} \label{convolution} For $f$ an $L^1$ function on $X$ and $g$ an $L^1$ function on $X$ which is radial
around the point $o$, the convolution of $f$ and $g$ is the function on $X$ defined by
$$
(f * g)(x) = \int_X f(y) (\tau_x g)(y) dvol(y)
$$
\end{definition}

\medskip

Note that, if $g = u \circ d_o$, then
\begin{align*}
||f * g ||_1 & \leq \int_X \int_X |f(y)| |(\tau_x g)(y)| dvol(y) dvol(x) \\
             & = \int_X |f(y)| \left(\int_X |u(d(x, y))| dvol(x)\right) dvol(y) \\
             & = \int_X |f(y)| \left(\int_{0}^{\infty} |u(r)| A(r) dr\right) dvol(y) \\
             & = ||f||_1 ||g||_1 \\
             & < +\infty \\
\end{align*}
so that the integral defining $(f * g)(x)$ exists for a.e. x, and $f * g \in L^1(X, dvol)$.

\medskip

Note if $f, g \in C^{\infty}_c(X)$ with $g = u \circ d_o$ radial around $o$, then $f * g$ is compactly supported. For the
Fourier transform of $f * g$ based at $o$, using the identity $B_{\xi,o}(x) = B_{\xi,o}(y) + B_{\xi, y}(x)$ we have
\begin{align*}
\widetilde{f * g}^o(\lambda, \xi) & = \int_X \left(\int_X f(y)u(d(x,y)) dvol(y)\right) e^{(-i\lambda - \rho)B_{\xi, o}(x)} dvol(x) \\
                                  & = \int_X f(y) e^{(-i\lambda - \rho)B_{\xi, o}(y)} \left(\int_X u(d(x,y)) e^{(-i\lambda - \rho)B_{\xi, y}(x)} dvol(x)\right) dvol(y) \\
                                  & = \int_X f(y) e^{(-i\lambda - \rho)B_{\xi, o}(y)} \widetilde{u \circ d_y}^y(\lambda, \xi) dvol(y) \\
                                  & = \tilde{f}^o(\lambda, \xi) \hat{u}(\lambda) \\
                                  & = \tilde{f}^o(\lambda, \xi) \hat{g}(\lambda) \\
\end{align*}
where we have used the fact that for the function $u \circ d_y$ which is radial around $y$ we have
$$
\widetilde{u \circ d_y}^y(\lambda, \xi) = \hat{u}(\lambda) = \hat{g}(\lambda)
$$
where $\hat{u}$ is the hypergroup Fourier transform of $u$ and $\hat{g}$ is the spherical Fourier transform of the function
$g$ which is radial around $o$.

\medskip

\begin{theorem} \label{radialconvolution} Let $L^1_o(X, dvol)$ denote the closed subspace of $L^1(X, dvol)$ consisting of those
$L^1$ functions which are radial around the point $o$. Then for $f, g \in L^1_o(X, dvol)$ we have $f * g \in L^1_o(X, dvol)$,
and $L^1_o(X, dvol)$ forms a commutative Banach algebra under convolution.
\end{theorem}

\medskip

\noindent{\bf Proof:} We first consider functions $f, g \in C^{\infty}_c(X)$ which are radial around $o$.
For $x \in X$, the function $M_x \phi_{\lambda, o}$ is an eigenfunction of $\Delta$ with eigenvalue $-(\lambda^2 + \rho^2)$
which is radial around the point
$x$ and takes the value $\phi_{\lambda, o}(x)$ at the point $x$, so it follows from Lemma \ref{uniqueeigenfn} that

$$
M_x \phi_{\lambda, o} = \phi_{\lambda, o}(x) \phi_{\lambda, x}
$$

The above equation together with the spherical Fourier inversion formula for the function $f$ gives, for $y \in Y$,
\begin{align*}
(M_x f)(y) & = C_0 \int_{0}^{\infty} (M_x \phi_{\lambda, o})(y) \hat{f}(\lambda) |c(\lambda)|^{-2} d\lambda \\
           & = C_0 \int_{0}^{\infty} \phi_{\lambda, o}(x) \phi_{\lambda, x}(y) \hat{f}(\lambda) |c(\lambda)|^{-2} d\lambda \\
\end{align*}

Now using the formal self-adjointness of $M_x$ we have
\begin{align*}
(f * g)(x) & = \int_X f(y) (\tau_x g)(y) dvol(y) \\
           & = \int_X f(y) (M_x \tau_x g)(y) dvol(y) \\
           & = \int_X (M_x f)(y) (\tau_x g)(y) dvol(y) \\
           & = C_0 \int_{0}^{\infty} \phi_{\lambda, o}(x) \hat{f}(\lambda) \left(\int_X  \phi_{\lambda, x}(y) (\tau_x g)(y) dvol(y)\right) |c(\lambda)|^{-2} d\lambda \\
           & = C_0 \int_{0}^{\infty} \phi_{\lambda, o}(x) \hat{f}(\lambda) \hat{g}(\lambda) |c(\lambda)|^{-2} d\lambda \\
\end{align*}

It follows from the above equation that $f * g$ is radial around the point $o$ since all the functions $\phi_{\lambda, o}$ are radial
around the point $o$, and it also follows that $f * g = g * f$.

\medskip

Now the inequality $||f * g||_1 \leq ||f||_1 ||g||_1$ implies, by the density of smooth, compactly supported radial functions in the space
$L^1_o(X, dvol)$, that for $f, g \in L^1_o(X, dvol)$ we have $f * g = g * f \in L^1_o(X, dvol)$, so $L^1_o(X, dvol)$ forms a
commutative Banach algebra under convolution. $\diamond$

\medskip

Finally, we remark that the radial hypergroup of the harmonic manifold $X$ can be realized as the convolution algebra
of finite radial measures on the manifold: convolution with radial measures can be defined, and the convolution of two radial
measures is again a radial measure. This can be proved by approximating finite radial measures by $L^1$ radial functions and
applying the previous Theorem. The convolution algebra $L^1_o(X, dvol)$ is then identified with a subalgebra of the hypergroup algebra
of finite radial measures under convolution.

\bibliography{moeb}
\bibliographystyle{alpha}

\end{document}